\pdfoutput=1
\documentclass[11pt]{amsart}

\usepackage[a4paper,margin=3.1cm]{geometry}
\usepackage[T1]{fontenc}
\usepackage{lmodern}
\usepackage{microtype}
\usepackage{mathtools,amssymb,amsthm}
\usepackage{enumitem}
\usepackage{needspace}
\usepackage{graphicx}
\usepackage{tikz}
\usetikzlibrary{arrows.meta,calc}
\usepackage{aliascnt}
\usepackage[hidelinks]{hyperref}
\usepackage[capitalise,nameinlink,noabbrev]{cleveref}
\hypersetup{
  pdftitle={Regular closedness off the real axis for root closures of bounded-height integer polynomials},
  pdfauthor={Bernat Espigul\'e},
  pdfsubject={Regular closedness of polynomial root closures and collinear connectedness loci},
  pdfkeywords={roots of bounded-height integer polynomials, regular closedness, connectedness locus, self-similar set, method of traps, finite capture}
}

\numberwithin{equation}{section}

\newtheorem{theorem}{Theorem}[section]
\newaliascnt{proposition}{theorem}
\newtheorem{proposition}[proposition]{Proposition}
\aliascntresetthe{proposition}
\newaliascnt{lemma}{theorem}
\newtheorem{lemma}[lemma]{Lemma}
\aliascntresetthe{lemma}
\newaliascnt{corollary}{theorem}
\newtheorem{corollary}[corollary]{Corollary}
\aliascntresetthe{corollary}
\theoremstyle{definition}
\newaliascnt{definition}{theorem}
\newtheorem{definition}[definition]{Definition}
\aliascntresetthe{definition}
\newaliascnt{example}{theorem}
\newtheorem{example}[example]{Example}
\aliascntresetthe{example}
\theoremstyle{remark}
\newaliascnt{remark}{theorem}
\newtheorem{remark}[remark]{Remark}
\aliascntresetthe{remark}

\crefname{theorem}{Theorem}{Theorems}
\Crefname{theorem}{Theorem}{Theorems}
\crefname{proposition}{Proposition}{Propositions}
\Crefname{proposition}{Proposition}{Propositions}
\crefname{lemma}{Lemma}{Lemmas}
\Crefname{lemma}{Lemma}{Lemmas}
\crefname{corollary}{Corollary}{Corollaries}
\Crefname{corollary}{Corollary}{Corollaries}
\crefname{definition}{Definition}{Definitions}
\Crefname{definition}{Definition}{Definitions}
\crefname{example}{Example}{Examples}
\Crefname{example}{Example}{Examples}
\crefname{remark}{Remark}{Remarks}
\Crefname{remark}{Remark}{Remarks}

\newcommand{\C}{\mathbb C}
\newcommand{\R}{\mathbb R}
\newcommand{\N}{\mathbb N_0}

\newcommand{\Mn}{\mathcal M_n}
\newcommand{\Core}{\mathcal K_n}
\newcommand{\An}{A_n}
\newcommand{\Bn}{B_n}
\newcommand{\Ktil}{\widetilde K}
\newcommand{\Rn}{\mathcal R_n}
\newcommand{\Xn}{\mathcal X_n}
\newcommand{\Thn}{\Theta_{\mathrm{fin}}(n)}
\newcommand{\Fill}[1]{\widehat{#1}}
\DeclareMathOperator{\conv}{conv}
\DeclareMathOperator{\dist}{dist}
\DeclareMathOperator{\diam}{diam}
\DeclareMathOperator{\interior}{int}
\DeclareMathOperator{\cl}{cl}
\DeclareMathOperator{\hop}{hop}
\renewcommand{\Re}{\operatorname{Re}}
\renewcommand{\Im}{\operatorname{Im}}
\newcommand{\doi}[1]{\href{https://doi.org/#1}{\nolinkurl{doi:#1}}}

\title[Regular closedness of polynomial root closures]
{Regular closedness off the real axis\\for root closures of bounded-height\\integer polynomials}

\author{Bernat Espigul\'e}
\address{Departament d'Inform\`atica, Matem\`atica Aplicada i Estad\'istica, Universitat de Girona, 17003 Girona, Catalonia, Spain}
\email{bernat@espigule.com}

\date{}

\begin{document}

\begin{abstract}
For each integer $h\ge1$, let $\mathcal Z_h$ be the closure in the complex plane of the roots of nonzero integer polynomials whose coefficients have absolute value at most $h$, with no degree restriction. We prove that $\mathcal Z_h$ is regular closed off the real axis: every non-real point of $\mathcal Z_h$ is a limit of interior points. The root-closure identification of Espigul\'e and Juher, for which we include a self-contained proof, shows that the same closure arises from monic polynomials and that outside the closed unit disc it is the connectedness locus of $h+1$ equally spaced collinear similarities. The main step proves that this locus is regular closed off the real axis; this was known for $h=1$ by Calegari, Koch and Walker, for $h\ge19$, and inside an explicit lens. The symmetry $z\mapsto1/z$ then gives the full-plane statement. The proof adapts the method of traps of Calegari, Koch and Walker, using a neighboring-piece criterion, short chains through consecutive digits, and Rouch\'e's theorem. We also classify the convex attractors.
\end{abstract}

\maketitle

\section{Introduction}\label{sec:introduction}

Fix an integer $h\ge1$.  The height of a nonzero polynomial $P(z)=\sum_{j=0}^m a_jz^j\in\mathbb Z[z]$ is $H(P):=\max_{0\le j\le m}|a_j|$.  Allowing arbitrary degrees, define the planar root closure
\begin{equation}\label{eq:planar-root-closure}
        \mathcal Z_h:=\cl_{\C}\{z\in\C:P(z)=0\text{ for some }0\ne P\in\mathbb Z[z]\text{ with }H(P)\le h\},
\end{equation}
where $\cl_{\C}$ and $\interior_{\C}$ denote closure and interior in the plane.  Our main result is the following.

\begin{theorem}[Regular closedness of polynomial root closures]\label{thm:root-closure}
For every integer $h\ge1$,
\begin{equation}\label{eq:global-regular-closedness}
        \mathcal Z_h\setminus\R
        =\cl_{\C}\bigl(\interior_{\C}(\mathcal Z_h)\bigr)\setminus\R.
\end{equation}
Equivalently, $\mathcal Z_h\setminus\R$ is regular closed in $\C\setminus\R$.
\end{theorem}

Thus every non-real point of the root closure can be approximated by interior points; the theorem concerns both sides of the unit circle.  Espigul\'e and Juher \cite[Proposition~1.1]{EJ2026FiniteCapture} identify $\mathcal Z_h$ with the closure obtained from monic polynomials and express it through the connectedness locus $\mathcal M_{h+1}$ of the collinear family introduced below and its image under $z\mapsto1/z$.  The new geometric step is regular closedness of $\mathcal M_{h+1}$ throughout its non-real part.  We state this step as \cref{thm:main} and deduce \cref{thm:root-closure} in \cref{sec:main-proof}, where \cref{prop:planar-identification} also gives a self-contained proof of the root-closure identification.  The restriction to non-real points is necessary: $0$ is an isolated point of $\mathcal Z_h$, and for $h\ge2$ the real points $z$ with $1+\sqrt h<|z|\le h+1$ lie in $\mathcal Z_h$ but are not limits of interior points, because near them $\mathcal Z_h$ coincides with $\mathcal M_{h+1}$ (\cref{prop:real-trace}\textup{(iii)}).

Connectedness loci, also called Mandelbrot sets, record the parameters for which an iterated function system has connected attractor.  For a pair of complex similarities, such a locus was introduced by Barnsley and Harrington \cite{BarnsleyHarrington1985}, and Bousch proved that it is connected and locally connected \cite{Bousch1992Thesis,Bousch1993}.  Bandt studied its topology and algorithms and asked whether every non-real connected parameter can be approximated by interior parameters \cite{Bandt2002}.  Solomyak and Xu proved this near the imaginary axis \cite{SolomyakXu2003}, and Calegari, Koch and Walker proved it in general by the method of traps \cite{CKW2017}.  Regular closedness is known for the connectedness loci of fractal $n$-gons with $n\ge3$: Bandt and Hung proved it for $n\ne4$ \cite{BandtHung2008} and Himeki and Ishii for $n=4$ \cite{HimekiIshii2020}; almost regular closedness for pairs of affine maps of the plane was proved by Rosler \cite{Rosler2025}.  Descriptions by zeros of power series were obtained by Solomyak and Xu \cite{SolomyakXu2003}, Shmerkin and Solomyak \cite{ShmerkinSolomyak2006}, and Solomyak \cite{Solomyak2015}.

We consider the following collinear family.  Fix an integer $n\ge2$ and let
\[
        \An:=\{0,1,\ldots,n-1\}.
\]
For a complex parameter $c$ with $|c|>1$, let $K(c,n)$ be the attractor of
\begin{equation}\label{eq:ifs}
        f_{j,c}(z)=j+\frac{z}{c},\qquad j\in\An,
\end{equation}
or equivalently
\begin{equation}\label{eq:series}
        K(c,n)=\left\{\sum_{k=0}^{\infty}a_kc^{-k}:a_k\in\An\right\}.
\end{equation}
For an infinite address $\alpha=(\alpha_0\alpha_1\ldots)\in\An^{\N}$, set
\[
        \pi_c(\alpha):=\sum_{j=0}^{\infty}\alpha_j c^{-j}.
\]
Here $\N=\{0,1,2,\ldots\}$, so addresses are indexed from zero.
Set
\[
        \Mn:=\{c\in\C:|c|>1\text{ and }K(c,n)\text{ is connected}\}.
\]
Closures, interiors and boundaries of subsets of $\{c\in\C:|c|>1\}$, such as $\Mn$, are taken relative to this open set, so relative interiors are planar interiors.  Operations on the planar root closure $\mathcal Z_h$ are explicitly marked by $\C$ and are taken in the full plane.

The digits $0,1,\ldots,n-1$ are convenient both geometrically and from the viewpoint of numeration.  They are affinely equivalent, for each fixed parameter, to the symmetric spacing-two alphabet used in \cite{EJS2024}; consequently the connectedness locus is unchanged.  The difference attractor is
\[
        D(c,n):=K(c,n)-K(c,n)
        =\left\{\sum_{k=0}^{\infty}d_k c^{-k}:d_k\in\{-(n-1),\ldots,n-1\}\right\}.
\]
This set contains the marked point $c$ precisely when $K(c,n)$ is connected.  In reciprocal coordinates this becomes a zero condition for a power series with bounded integer coefficients.  \Cref{prop:neighbor} proves this marked-point criterion; its proof includes a planar crossing argument showing that connectedness forces overlap of neighboring first-level pieces.

For the collinear family, Espigul\'e, Juher and Salda\~na \cite[Theorem~2]{EJS2024} proved that every non-real connected parameter in the closed lens $|c\pm1|\le\sqrt{2n}$ is a limit of interior parameters, and established the same conclusion throughout the non-real locus when $n\ge21$.  They conjectured regular closedness of the off-lens part $\Mn\setminus(\R\cup\overline{\Xn})$ for $3\le n\le20$ \cite[Section~5]{EJS2024}; here $\overline{\Xn}$ is the closed lens, the closure in parameter space of the open lens $\Xn$ defined in \eqref{eq:lens} below.  Espigul\'e and Juher added the case $n=20$, with an explicit dense set of certified interior parameters formed by semialgebraic capture levels \cite[Corollary~5.4]{EJ2026FiniteCapture}, and showed that the lens contains the whole non-real locus exactly when $n\ge20$ \cite[Theorem~5.3]{EJ2026FiniteCapture}.  The theorem below proves regular closedness away from the real axis for every $n\ge2$, without a lens restriction; for $n=2$ it is the theorem of \cite{CKW2017}, for $n\ge20$ it follows from \cite{EJS2024,EJ2026FiniteCapture}, and the cases $3\le n\le19$ are new.

The proof uses traps in the sense of Calegari, Koch and Walker \cite{CKW2017}: stable crossings in enclosures of cylinder sets force overlap. These topological traps differ from the canonical trap of \cite{EJ2026FiniteCapture}, a self-covering parallelogram that contains the covering rectangle of \cite{EJS2024}: entry of a finite inverse orbit of the marked point into the canonical trap certifies an interior parameter, whereas a topological trap forces overlap without necessarily providing such an entry region. The adaptation uses a neighboring-piece criterion, short chains through consecutive digits, and a classification of convex attractors.

\begin{theorem}\label{thm:main}
For every integer $n\ge2$,
\begin{equation}\label{eq:main}
        \Mn\setminus\R
        =\overline{\interior(\Mn)}\setminus\R.
\end{equation}
Equivalently, $\Mn\setminus\R$ is regular closed in $\{c\in\C\setminus\R:|c|>1\}$.
\end{theorem}

\begin{figure}[t]
\centering
\includegraphics[width=\textwidth]{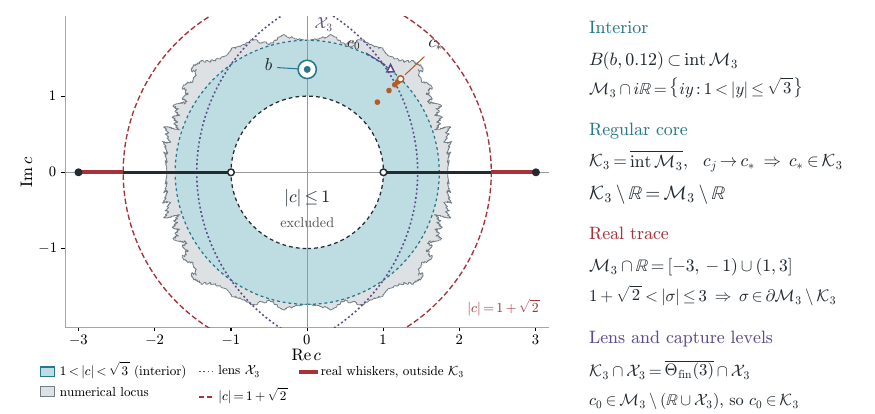}
\caption{The connectedness locus $\mathcal M_3$ and its regular core.
Blue: the annulus $1<|c|<\sqrt3$, interior by \cref{lem:annulus} and optimal
on the imaginary axis by \cref{prop:real-trace}\textup{(iv)}, with the
interior disc $B(b,0.12)$, $b=1.35i$, and a radial sequence converging to
$c_*=\sqrt3e^{\pi i/4}\in\mathcal K_3$.
Black and red: the real trace $[-3,-1)\cup(1,3]$ of \cref{prop:real-trace};
its red part $1+\sqrt2<|\sigma|\le3$ lies outside the core, since
$\mathcal M_3\setminus\R$ is inside the red circle $|c|=1+\sqrt2$.
Dotted: the boundary of the lens $\mathcal X_3$, inside which the core is the closure
of the capture levels $\Theta_k(3)$ (\cref{cor:finite-capture}).
The triangle is the parameter $c_0$ of \cref{fig:finite-trap} and
\cref{ex:worked}, a polynomial root outside both the lens and the annulus.}
\label{fig:regular-core}
\end{figure}

An interior parameter has a neighborhood throughout which the attractor remains connected. Thus \cref{thm:main} says that every non-real connected parameter can be approximated by parameters with this stability property. We call
\[
        \Core:=\overline{\interior(\Mn)}
\]
the \emph{regular core}. Any point of $\Mn$ outside $\Core$ lies on the real axis. In \cref{cor:regular-core} we identify $\Core$ with the closure of the non-real roots of modulus greater than $1$ of monic integer polynomials of height at most $n-1$, and we show that the non-real boundary of $\Mn$ is the boundary of its interior.  In \cref{sec:structure} we describe the internal structure of the locus that results: the real and imaginary traces, the real parameters that leave the core, and, inside the lens, the identification of the core with the closure of the semialgebraic capture levels of \cite{EJ2026FiniteCapture}. Figure~\ref{fig:regular-core} illustrates these statements for $n=3$.

We also determine the convexity locus.  Convexity occurs exactly on rational-direction radial segments: if $c=\rho e^{\pi i p/q}$ with $p/q$ reduced, then the attractor is convex if and only if $\rho^q\le n$.  Convexity does not occur in irrational directions.  For $n=2$, this specializes to \cite[Lemma~7.2.3]{CKW2017}.  For general $n$, the necessity of $\rho^q\le n$ and the absence of convexity in irrational directions also appear in Silvestri and P\'erez \cite[Section~3.5]{SilvestriPerez2026}, which describes the extreme points of planar self-affine attractors with a common linear part by sofic systems; \cref{thm:convex-locus} adds the converse and the equivalence with convexity of the filled hull (defined below), which is the form the trap argument uses.

For a compact set $L\subset\C$, its \emph{filled hull} $\Fill L$ is the union of $L$ and all bounded components of $\C\setminus L$. A \emph{continuum} is a nonempty compact connected set; it is \emph{full} if its complement is connected.

The proof of \cref{thm:main} is organized as follows.  The annulus $1<|c|<\sqrt n$ lies in $\interior(\Mn)$ by a similarity-dimension argument (\cref{lem:annulus}).  For $|c|>\sqrt n$, the convexity theorem of \cref{sec:convex} implies that a non-real connected attractor has nonconvex filled hull.  The planar trap lemma (\cref{prop:stable-linked}, proved in \cref{app:linked-translates}) then provides a stable crossing configuration for a translate of that hull.  A short-hop argument converts this configuration into overlap of adjacent cylinders (\cref{sec:traps}), and Rouch\'e's theorem realizes the required cylinder translation at parameters arbitrarily close to the original one (\cref{lem:rouche}).  \Cref{sec:main-proof} completes the proofs of \cref{thm:main,thm:root-closure}, and \cref{sec:structure} combines them with the closure of polynomial roots and with the capture levels of \cite{EJ2026FiniteCapture} to describe the internal structure of $\Mn$.

\section{Neighbor overlaps and adjacent translations}\label{sec:prelim}

\subsection{Affine normalization and the marked point}

We first record the affine equivalence with the symmetric alphabet used in \cite{EJS2024}.

\begin{lemma}[Affine normalization]\label{lem:normalization}
Let
\[
        A_n^{\mathrm{sym}}:=\{-n+1,-n+3,\ldots,n-1\}
\]
and let $E(c,n)$ be the corresponding attractor.  Then
\begin{equation}\label{eq:affine-normalization}
        E(c,n)=2K(c,n)-(n-1)\frac{c}{c-1}.
\end{equation}
In particular, $E(c,n)$ and $K(c,n)$ are affinely homeomorphic for each $c$ with $|c|>1$ and have the same connectedness locus.
\end{lemma}

\begin{proof}
Write each symmetric digit as $2a-(n-1)$ with $a\in\An$ and sum the resulting geometric series.
\end{proof}

We use the following standard planar crossing fact; the horizontal-chord lemma derived from it is the additive part of Hopf's theorem on chords of plane continua~\cite[Satz~I]{Hopf1936}.

\begin{lemma}[Rectangle crossing]\label{lem:rectangle-crossing}
Let $R=[a,b]\times[u,v]$ be a nondegenerate closed rectangle.  If $H,V\subset R$ are compact connected sets such that $H$ meets both vertical sides of $R$ and $V$ meets both horizontal sides, then $H\cap V\ne\varnothing$.
\end{lemma}

\begin{proof}
For paths in place of continua this is \cite[Lemma~2]{Maehara1984}, the lemma underlying Maehara's proof of the Jordan curve theorem, where it is derived from Brouwer's fixed-point theorem; the general case reduces to it.  Suppose that $H\cap V=\varnothing$ and let $\varepsilon:=\dist(H,V)>0$.  Subdivide $R$ into finitely many closed rectangles of diameter less than $\varepsilon/2$, and let $H'$ be the union of those that meet $H$ and $V'$ the union of those that meet $V$.  If a rectangle $A$ meeting $H$ intersected a rectangle $B$ meeting $V$, points $h\in A\cap H$ and $v\in B\cap V$ would satisfy $|h-v|\le\diam(A)+\diam(B)<\varepsilon$, a contradiction.  Hence $H'\cap V'=\varnothing$.  The set $H'$ is connected, being the union of the connected set $H$ with closed rectangles each of which meets $H$; as a connected finite union of closed rectangles it is path connected, and it meets both vertical sides of $R$ because it contains $H$.  Hence $H'$ contains a path joining the two vertical sides of $R$, and likewise $V'$ contains a path joining the two horizontal sides.  These paths are disjoint, contradicting the case of paths.  The lemma is also a consequence of the Phragm\'en--Brouwer separation property; see, for example, \cite[Lemma~7]{PascolettiZanolin2013}.
\end{proof}

\Needspace{10\baselineskip}
\begin{lemma}[Horizontal chords of a plane continuum]\label{lem:horizontal-chords}
Let $L\subset\C$ be nonempty, compact and connected, and let $a,b>0$.  If
\[
        L\cap\bigl(L+(a+b)\bigr)\ne\varnothing,
\]
then
\[
        L\cap(L+a)\ne\varnothing
        \quad\text{or}\quad
        L\cap(L+b)\ne\varnothing.
\]
Equivalently, calling $a>0$ a \emph{missing chord length} of $L$ when $L\cap(L+a)=\varnothing$, the set of missing chord lengths is closed under addition.
\end{lemma}

\begin{proof}
Let
\[
 x_-:=\min_{z\in L}\Re z,\quad x_+:=\max_{z\in L}\Re z,
 \qquad
 y_-:=\min_{z\in L}\Im z,\quad y_+:=\max_{z\in L}\Im z.
\]
If $y_-=y_+$, then $L$ is a compact connected subset of a horizontal line and hence an interval, so the assertion follows.  Assume $y_-<y_+$ and put
\[
        R=[x_-,x_++a+b]\times[y_-,y_+].
\]
The set
\[
        H:=L\cup\bigl(L+(a+b)\bigr)
\]
is a continuum because its two terms intersect; it lies in $R$ and meets both vertical sides.  The continuum $V:=L+a$ lies in $R$ and meets both horizontal sides.  By \Cref{lem:rectangle-crossing}, $H\cap V\ne\varnothing$.  Intersection with $L$ gives $L\cap(L+a)\ne\varnothing$; intersection with $L+(a+b)$, followed by translation by $-a$, gives $L\cap(L+b)\ne\varnothing$.
\end{proof}

Let $S(c):=c^{-1}K(c,n)$.  The first-level pieces are $j+S(c)$, $j\in\An$, and we set
\begin{equation}\label{eq:delta}
        \delta_n(c):=\dist\bigl(S(c),S(c)+1\bigr).
\end{equation}
By translation, this is the distance between any neighboring pair of first-level pieces.

The next lemma is the collinear form of the short hop lemma of Calegari, Koch and Walker \cite[Lemma~5.2.2]{CKW2017}; the hop length is controlled by the neighboring-piece distance $\delta_n(c)$.

\begin{lemma}[Collinear short-hop chains]\label{lem:short-hop}
Let $c\in\C$ with $|c|>1$, put $K:=K(c,n)$ and $\delta:=\delta_n(c)$.  For every $\tau>0$ and every $x,y\in K$, there is a polygonal path from $x$ to $y$ whose vertices lie in $K$ and whose edges have length less than $\delta+\tau$.
\end{lemma}

\begin{proof}
For $r>0$, an $r$-chain is a finite sequence of points with consecutive distances less than $r$.  For $x,y\in K$ let $\hop(x,y)$ be the infimum of the numbers $r>0$ for which an $r$-chain in $K$ joins $x$ to $y$, and put
\[
        \hop_K:=\sup_{x,y\in K}\hop(x,y)\le\diam K<\infty.
\]
Fix $\tau'>0$.  If $x$ and $y$ lie in the same first-level piece $f_{j,c}(K)$, write $x=f_{j,c}(x')$ and $y=f_{j,c}(y')$ with $x',y'\in K$; the image under $f_{j,c}$ of a $(\hop_K+\tau')$-chain in $K$ from $x'$ to $y'$ is an $|c|^{-1}(\hop_K+\tau')$-chain in $f_{j,c}(K)$ from $x$ to $y$.  If $x\in f_{j,c}(K)$ and $y\in f_{k,c}(K)$ with $j<k$, after interchanging $x$ and $y$ if necessary, choose for each $i\in\{j,\ldots,k-1\}$ points $z_i\in f_{i,c}(K)$ and $z_i'\in f_{i+1,c}(K)$ with $|z_i-z_i'|=\delta$, which exist by compactness, and concatenate chains of the previous kind inside the pieces $f_{j,c}(K),\ldots,f_{k,c}(K)$, from $x$ to $z_j$, from $z_i'$ to $z_{i+1}$, and from $z_{k-1}'$ to $y$, with the jumps from $z_i$ to $z_i'$.  Every jump of the resulting chain is less than $\max\{\delta,|c|^{-1}\hop_K\}+\tau'$.  Hence $\hop_K\le\max\{\delta,|c|^{-1}\hop_K\}+\tau'$ for every $\tau'>0$, so
\[
        \hop_K\le\max\{\delta,|c|^{-1}\hop_K\},
\]
and since $|c|^{-1}<1$ this forces $\hop_K\le\delta$.  Consequently any two points of $K$ are joined by a $(\delta+\tau)$-chain in $K$; joining its consecutive points by line segments gives the required polygonal path.
\end{proof}

The equivalence of \textup{(i)}--\textup{(iii)} below, the neighboring-overlap criterion for collinear digits, is due to Bousch \cite[Chapter~1, p.~41]{Bousch1992Thesis}; we include a proof because the same chord argument also gives the final assertion.

\begin{proposition}[Neighboring-piece and marked-point criterion]\label{prop:neighbor}
For $c\in\C$ with $|c|>1$, the following are equivalent:
\begin{enumerate}[label=\textup{(\roman*)},leftmargin=2.2em]
\item $K(c,n)$ is connected;
\item $S(c)\cap\bigl(S(c)+1\bigr)\ne\varnothing$;
\item every neighboring pair of first-level pieces intersects;
\item $c\in D(c,n)$;
\item there exist addresses $\alpha,\beta\in\An^{\N}$ with $|\alpha_0-\beta_0|=1$ and $\pi_c(\alpha)=\pi_c(\beta)$.
\end{enumerate}
If these conditions fail, all first-level pieces are pairwise disjoint.
\end{proposition}

\begin{proof}
The equivalence of \textup{(ii)} and \textup{(iii)} follows by translation.  If \textup{(ii)} holds, the first-level intersection graph contains the path through the ordered digits $0,1,\ldots,n-1$; Hata's criterion therefore gives \textup{(i)} \cite{Hata1985}.

Conversely, suppose that $K(c,n)$ is connected.  Then $S(c)$ is a plane continuum.  If the chord length $1$ were missing from $S(c)$, the additive formulation of \cref{lem:horizontal-chords} would imply inductively that every positive integer chord length is missing.  Thus the translates $j+S(c)$, $j\in\An$, would be pairwise disjoint.  Their finite union is $K(c,n)$, contradicting connectedness.  Hence \textup{(ii)} holds.

Condition \textup{(ii)} means that for some $x,y\in K(c,n)$,
\[
        \frac{x}{c}=1+\frac{y}{c},
\]
which is equivalent to $c=x-y\in D(c,n)$.  This proves the equivalence with \textup{(iv)}.  If $x=\pi_c(\alpha')$ and $y=\pi_c(\beta')$, then the addresses $\alpha:=0\alpha'$ and $\beta:=1\beta'$ satisfy $\pi_c(\alpha)=x/c=1+y/c=\pi_c(\beta)$, which gives \textup{(v)}; reversing the argument proves the converse.

To prove the final assertion, suppose that these conditions fail and put
$\delta=\delta_n(c)>0$.  Choose
\[
        \frac{\delta}{2|c|}<r<\frac{\delta}{2}
\]
and let $L:=\{z\in\C:\dist(z,S(c))\le r\}$.  Take $\tau>0$ so small that
$(\delta+\tau)/|c|<2r$.  By \cref{lem:short-hop} applied to $K(c,n)$ and scaled by $c^{-1}$, any two points of $S(c)$
can be joined by a chain in $S(c)$ with jumps less than $2r$.  The closed
discs of radius $r$ centered along such a chain overlap consecutively.
Their union over all centers in $S(c)$ is $L$, so $L$ is a compact
continuum.  Moreover,
\[
        \dist(L,L+1)\ge\delta-2r>0.
\]
By \cref{lem:horizontal-chords}, every positive integer chord length is
missing from $L$, and hence from $S(c)$.  Thus the first-level pieces
$j+S(c)$ are pairwise disjoint.
\end{proof}

The last assertion is a dichotomy particular to collinear digit sets: the first-level intersection graph of $K(c,n)$ either contains the path through the consecutive digits $0,1,\ldots,n-1$ or has no edges at all.  Hata's criterion by itself would permit connected intersection graphs in which some neighboring pairs of pieces are disjoint; \cref{lem:horizontal-chords} excludes them.  In particular $K(c,n)$ is either connected or totally disconnected: when the first-level pieces are pairwise disjoint, so are the level-$m$ pieces $\sum_{j<m}u_jc^{-j}+c^{-m}K(c,n)$, $u\in\An^m$, for every $m$, and their diameters tend to zero.  This is what makes the single quantity $\delta_n(c)$, which vanishes exactly on $\Mn$, the right measure of disconnectedness for the short-hop and trap arguments of \cref{sec:traps}.

\begin{corollary}[Power-series formulation]\label{cor:power-series}
A parameter $c$ with $|c|>1$ belongs to $\Mn$ if and only if there are coefficients
\[
        d_k\in\{-(n-1),\ldots,n-1\}
\]
such that
\begin{equation}\label{eq:power-series}
        1-\sum_{k=0}^{\infty}d_kc^{-(k+1)}=0.
\end{equation}
In particular, $\Mn$ is closed in $\{z\in\C:|z|>1\}$.
\end{corollary}

\begin{proof}
The equivalence follows from $c\in D(c,n)$ after division by $c$.  The attractor $D(c,n)$ depends continuously on $c$ in the Hausdorff metric by Hutchinson's theorem \cite{Hutchinson1981}, the maps of the system being uniformly contractive near each parameter; hence the relation $c\in D(c,n)$ is closed.
\end{proof}

\begin{lemma}[Interior annulus]\label{lem:annulus}
For every $n\ge2$,
\begin{equation}\label{eq:annulus}
        \{c\in\C:1<|c|<\sqrt n\}\subset\interior(\Mn).
\end{equation}
\end{lemma}

\begin{proof}
Suppose $1<|c|<\sqrt n$ and $K(c,n)$ is disconnected.  By \Cref{prop:neighbor}, the first-level pieces are pairwise disjoint, so the IFS satisfies the strong separation condition.  Its Hausdorff dimension therefore equals its similarity dimension \cite{Hutchinson1981}:
\[
        \dim_{\mathrm H}K(c,n)=\frac{\log n}{\log|c|}>2,
\]
a contradiction.  Thus the open annulus in \eqref{eq:annulus} is contained in $\Mn$ and hence in its interior.
\end{proof}

The annulus was obtained in \cite[Proposition~5(i)]{EJS2024}; for $n=2$ it is the inner annulus of \cite[Lemma~3.2.1]{CKW2017}.  The radius $\sqrt n$ is optimal: by \cref{prop:real-trace}\textup{(iv)} below, the locus meets the imaginary axis exactly in the points $iy$ with $1<|y|\le\sqrt n$.  In the reciprocal variable $\lambda=1/c$, \cref{lem:annulus}, this optimality and \cref{cor:power-series} show that $n^{-1/2}\le|\lambda|<1$ is the largest annulus of the form $r\le|\lambda|<1$ consisting of zeros of power series $1+\sum_{k\ge1}e_k\lambda^k$ with $e_k\in\{-(n-1),\ldots,n-1\}$.  Kewarth \cite{Kewarth2026} proved that the smaller annulus $(n-1)^{-1/2}\le|\lambda|<1$ has this property, noted the bound on the imaginary axis, and suggested that $n^{-1/2}$ is the optimal radius.

\subsection{Adjacent cylinder translations}

For an address $\alpha$ and $m\ge1$, write $\alpha|_m:=\alpha_0\cdots\alpha_{m-1}$ for its prefix of length $m$.  For a finite word $u=u_0\cdots u_{m-1}\in\An^m$, set
\[
        \pi_c(u):=\sum_{j=0}^{m-1}u_j c^{-j},
        \qquad
        K_u(c):=\pi_c(u)+c^{-m}K(c,n).
\]

\begin{definition}[Normalized cylinder translation]\label{def:translation}
For words $u,v\in\An^m$ of the same length, define
\begin{equation}\label{eq:translation}
        w_{u,v}(c):=c^m\bigl(\pi_c(u)-\pi_c(v)\bigr).
\end{equation}
The pair $(u,v)$ is \emph{adjacent} if $|u_0-v_0|=1$.
\end{definition}

\begin{lemma}[Cylinder translation identities]\label{lem:cylinder-identities}
Let $u,v\in\An^m$ and put $w=w_{u,v}(c)$.  Then
\[
        K_u(c)\cap K_v(c)\ne\varnothing
        \quad\Longleftrightarrow\quad
        K(c,n)\cap\bigl(K(c,n)+w\bigr)\ne\varnothing,
\]
and
\begin{equation}\label{eq:distance-identity}
        \dist\bigl(K(c,n),K(c,n)+w\bigr)
        =|c|^m\dist\bigl(K_u(c),K_v(c)\bigr).
\end{equation}
\end{lemma}

\begin{proof}
Write $K=K(c,n)$.  If $z\in K$ and $z-w\in K$, then
\[
        \pi_c(v)+c^{-m}z=\pi_c(u)+c^{-m}(z-w),
\]
which gives a point of $K_v(c)\cap K_u(c)$.  The converse is the same calculation in reverse.  For \eqref{eq:distance-identity}, multiply the cylinder distance by $|c|^m$ and use the symmetry of the difference set $K-K$.
\end{proof}

\begin{lemma}[Rouch\'e realization of adjacent translations]\label{lem:rouche}
Let $c_0\in\Mn$ and choose $\alpha,\beta\in\An^{\N}$, as in \cref{prop:neighbor}\textup{(v)}, with
\[
        |\alpha_0-\beta_0|=1,
        \qquad
        \pi_{c_0}(\alpha)=\pi_{c_0}(\beta).
\]
For every $w\in\C$ and every neighborhood $U$ of $c_0$ in $\{z\in\C:|z|>1\}$, there exist $m\ge1$ and $c\in U$ such that
\[
        w_{\alpha|_m,\beta|_m}(c)=w.
\]
\end{lemma}

\begin{proof}
Define holomorphic functions on $\{z\in\C:|z|>1\}$ by
\[
        T(c):=\pi_c(\alpha)-\pi_c(\beta),
        \qquad
        T_m(c):=\pi_c(\alpha|_m)-\pi_c(\beta|_m).
\]
The function $T$ is not identically zero because $T(c)\to\alpha_0-\beta_0\ne0$ as $|c|\to\infty$, while $T(c_0)=0$.  Moreover $T_m\to T$ uniformly on compact subsets of $\{z\in\C:|z|>1\}$.

Choose a closed disc $V\Subset U$ centered at $c_0$ such that $T$ has no zero on $\partial V$.  Put
\[
        H_m(c):=T_m(c)-wc^{-m}.
\]
Both $T_m-T$ and $wc^{-m}$ tend uniformly to zero on $V$.  Rouch\'e's theorem therefore implies, for all sufficiently large $m$, that $H_m$ and $T$ have the same positive number of zeros in the interior of $V$, counted with multiplicity.  At a zero $c$ of $H_m$,
\[
        c^m T_m(c)=w,
\]
which is the required identity.  The two prefixes remain adjacent because their first letters are $\alpha_0$ and $\beta_0$.
\end{proof}

\Cref{lem:rouche} is the collinear counterpart of the surjective perturbation lemma of \cite[Lemma~7.2.5 and Corollary~7.2.6]{CKW2017}.

\section{The convexity locus}\label{sec:convex}

Translate the attractor to its centered form
\begin{equation}\label{eq:centered}
        \Ktil(c,n):=K(c,n)-\frac{n-1}{2}\frac{c}{c-1}
        =\left\{\sum_{k=0}^{\infty}b_kc^{-k}:b_k\in\Bn\right\},
\end{equation}
where
\[
        \Bn:=\left\{j-\frac{n-1}{2}:j=0,1,\ldots,n-1\right\},
        \qquad
        J_n:=\left[-\frac{n-1}{2},\frac{n-1}{2}\right].
\]
Convexity is unchanged by this translation.

\begin{lemma}[Boundary of the filled hull]\label{lem:filled-boundary}
For every compact $L\subset\C$,
\[
        \partial\Fill{L}\subset L.
\]
\end{lemma}

\begin{proof}
If $U_\infty$ is the unbounded component of $\C\setminus L$, then $\Fill L=\C\setminus U_\infty$ and $\partial\Fill L=\partial U_\infty\subset L$.
\end{proof}

\begin{lemma}[Convex-hull formula]\label{lem:convex-hull}
For every $c\in\C$ with $|c|>1$,
\begin{equation}\label{eq:minkowski}
        \conv\Ktil(c,n)=J_n+c^{-1}\conv\Ktil(c,n)=\sum_{j=0}^{\infty}c^{-j}J_n,
\end{equation}
where the last expression is a Hausdorff-convergent Minkowski sum of line segments.
\end{lemma}

\begin{proof}
The self-similarity equation $\Ktil(c,n)=\Bn+c^{-1}\Ktil(c,n)$ gives the first identity after taking convex hulls, since the convex hull of a Minkowski sum is the Minkowski sum of the convex hulls and $\conv\Bn=J_n$.  Iteration and $|c|>1$ yield the second.
\end{proof}

\begin{lemma}[Real interval criterion]\label{lem:real-criterion}
Let $\sigma\in\R$ with $|\sigma|>1$.  Then $\Ktil(\sigma,n)$ is an interval if and only if $|\sigma|\le n$.  If $|\sigma|>n$, the attractor contains no nondegenerate interval.
\end{lemma}

\begin{proof}
Because $\Bn$ is symmetric, replacing $b_k$ by $(-1)^k b_k$ shows that $\Ktil(-\rho,n)=\Ktil(\rho,n)$ for $\rho>1$.  It therefore suffices to take $\sigma=\rho>1$.

The convex hull is $[-M,M]$, where
\[
        M=\frac{(n-1)\rho}{2(\rho-1)}.
\]
Its first-level images have centers separated by $1$ and length $2M/\rho=(n-1)/(\rho-1)$.  They cover $[-M,M]$ without gaps exactly when $\rho\le n$, in which case the interval is invariant and hence is the attractor.

If $\rho>n$, the first-level convex-hull intervals have positive gaps, and the same is true inside every cylinder.  At level $k$, the attractor is contained in pairwise disjoint intervals of length $2M\rho^{-k}$.  A nondegenerate interval contained in the attractor would have to lie in one level-$k$ interval for every $k$, which is impossible as these lengths tend to zero.
\end{proof}

\begin{theorem}[Convexity locus]\label{thm:convex-locus}
Let $c=\rho e^{i\theta}\in\C$ with $\rho>1$.  The following are equivalent:
\begin{enumerate}[label=\textup{(\roman*)},leftmargin=2.2em]
\item $K(c,n)$ is convex;
\item the filled hull $\Fill{K(c,n)}$ is convex;
\item $\theta/\pi$ is rational and, writing $\theta/\pi=p/q$ in lowest terms with $q\ge1$, one has $\rho^q\le n$.
\end{enumerate}
Consequently,
\begin{equation}\label{eq:convex-locus}
 \{c\in\C:|c|>1\text{ and }K(c,n)\text{ is convex}\}
 =\bigcup_{\substack{p\in\mathbb Z,\ q\ge1\\ \gcd(p,q)=1}}
   \{\rho e^{\pi i p/q}:1<\rho\le n^{1/q}\}.
\end{equation}
\end{theorem}

\begin{proof}
A convex compact set equals its filled hull, so \textup{(i)} implies \textup{(ii)}.  Assume \textup{(ii)} and work with the centered attractor $L:=\Ktil(c,n)$.  Put $X:=\Fill L$.  Since $X$ is convex and contains $L$, $\conv L\subset X$.  Conversely, $\Fill L\subset\conv L$ for every compact planar set $L$, so $X=\conv L$.

If $c\in\R$, then $L\subset\R$ and $\Fill L=L$, so the assertion follows from \cref{lem:real-criterion}, with $q=1$.  Assume $c\notin\R$.  By \Cref{lem:convex-hull}, $X=J_n+c^{-1}X$.  For compact convex sets $A,B\subset\C$ and a linear functional $\phi$, the face of $A+B$ on which $\phi$ is maximal is the Minkowski sum of the faces of $A$ and of $B$ on which $\phi$ is maximal: indeed, $\phi(a+b)=\phi(a)+\phi(b)\le\max_A\phi+\max_B\phi$, with equality exactly when $a$ and $b$ attain the two maxima.  Applied with $\phi=\Im$, $A=J_n$ and $B=c^{-1}X$, and since $\Im$ vanishes identically on the horizontal segment $J_n$, this shows that the upper support face of $X$ for $\Im$ contains a translate of the nondegenerate segment $J_n$.  By \Cref{lem:filled-boundary}, this support segment lies in $L$.

If $\theta/\pi$ were irrational, then $\Im(c^{-j})\ne0$ for every $j\ge1$.  Maximizing the imaginary part of
\[
        \sum_{j=0}^{\infty}b_j c^{-j},\qquad b_j\in\Bn,
\]
forces a unique endpoint choice for every $b_j$ with $j\ge1$; only the $n$ possible values of $b_0$ remain.  Thus $L$ would meet its upper support line in at most $n$ points, contradicting the segment above.  Hence $\theta/\pi=p/q\in\mathbb Q$ in lowest terms.

The indices at which $c^{-j}$ is real are exactly the multiples of $q$.  On the upper support line, all other digits are forced, while the digits with indices divisible by $q$ are free.  The support slice of $L$ is therefore a translate of
\[
        \Ktil(c^q,n)=\Ktil\bigl((-1)^p\rho^q,n\bigr).
\]
Since the slice contains a nondegenerate interval, \cref{lem:real-criterion} gives $\rho^q\le n$.  This proves \textup{(ii)}$\Rightarrow$\textup{(iii)}.

Finally assume \textup{(iii)}.  Group the defining series for $\Ktil(c,n)$ by residues modulo $q$:
\[
        \Ktil(c,n)=\sum_{r=0}^{q-1}c^{-r}\Ktil(c^q,n).
\]
The set $\Ktil(c^q,n)$ is an interval by \cref{lem:real-criterion}.  Hence $\Ktil(c,n)$ is a finite Minkowski sum of line segments and is convex.  This proves \textup{(iii)}$\Rightarrow$\textup{(i)} and \eqref{eq:convex-locus}.
\end{proof}

\begin{corollary}[Finite polygonal representation]\label{cor:finite-polygon}
Suppose $c=\rho e^{\pi i p/q}$, where $p/q$ is reduced and $\rho^q\le n$.  Then
\begin{equation}\label{eq:finite-polygon}
        \Ktil(c,n)
        =\frac{1}{1-\rho^{-q}}
          \sum_{r=0}^{q-1}c^{-r}J_n.
\end{equation}
In particular, $\Ktil(c,n)$ is a centrally symmetric convex polygon, possibly degenerate to a line segment.
\end{corollary}

\begin{proof}
Group the defining series by residues modulo $q$.  Since $c^q=(-1)^p\rho^q$ with $\rho^q\le n$, \cref{lem:real-criterion} shows that $\Ktil(c^q,n)$ is an interval, hence equals its convex hull, and, since $J_n=-J_n$, \cref{lem:convex-hull} gives
\[
        \Ktil(c^q,n)
        =\sum_{k=0}^{\infty}c^{-qk}J_n
        =\frac{1}{1-\rho^{-q}}J_n.
\]
Substitution into
$\Ktil(c,n)=\sum_{r=0}^{q-1}c^{-r}\Ktil(c^q,n)$
yields \eqref{eq:finite-polygon}.
\end{proof}

\begin{figure}[b]
\centering
\input{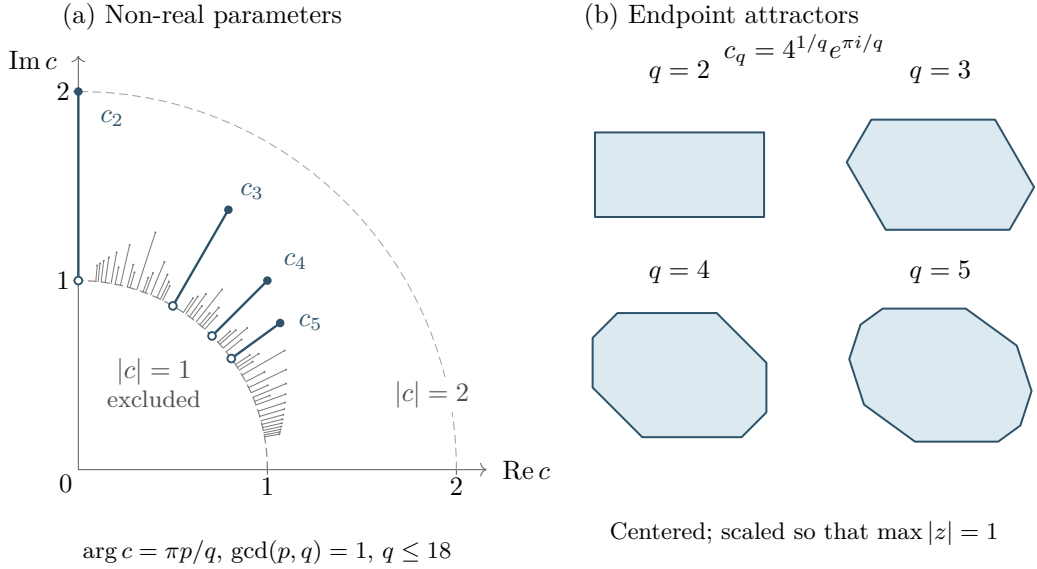}
\caption{The convexity locus for $n=4$.
\textup{(a)} The radial segments $1<|c|\le4^{1/q}$ in every reduced
direction $\arg c=\pi p/q$ with $0<p/q\le1/2$ and $q\le18$,
as given by \cref{thm:convex-locus}; real parameters are omitted.
The unit circle is excluded.
The endpoints $c_q=4^{1/q}e^{\pi i/q}$ for $q=2,3,4,5$ are marked.
\textup{(b)} The corresponding centered attractors $\Ktil(c_q,4)$,
computed from \eqref{eq:finite-polygon} and scaled so that
$\max |z|=1$, with their orientations preserved.}
\label{fig:convexity-locus}
\end{figure}

\Cref{fig:convexity-locus} shows the convexity locus for $n=4$ and the polygons \eqref{eq:finite-polygon} at the endpoints of four of its segments.

\begin{corollary}[Convex filled hulls off the real axis]\label{cor:convex-bound}
Let $c\in\C\setminus\R$ with $|c|>1$.  If $\Fill{K(c,n)}$ is convex, then
\[
        |c|\le\sqrt n.
\]
\end{corollary}

\begin{proof}
In \Cref{thm:convex-locus}, a non-real direction has $q\ge2$, so $|c|\le n^{1/q}\le\sqrt n$.
\end{proof}

\section{Linked translations and short-hop stability}\label{sec:traps}

We formulate the crossing and proximity conditions used in the trap method of \cite[Section~7.1]{CKW2017}. After applying a common similarity to two equal-depth cylinder sets, their enclosures become translates of one disc. The configurations constructed below have nonzero algebraic intersection number, as in their definition; the connectedness argument needs only the crossing property stated here.

For $r>0$ and $A\subset\C$, write
\[
        \mathcal N_r(A):=\{z\in\C:\dist(z,A)\le r\}.
\]
The following proposition supplies the topological configuration used in the method of traps; it is a stable form of the trap-like vector of \cite[Definition~7.2.1 and Lemma~7.2.2]{CKW2017}, with the neighborhoods and the perturbation range made explicit.  Its proof is given in \cref{app:linked-translates}.

\begin{proposition}[Stable linked translate]\label{prop:stable-linked}
Let $X\subset\C$ be a nonconvex full continuum.  There exist a vector $w\in\C$, a compact topological disc $Q$, numbers $\varepsilon_0,\eta_0>0$, points
\[
        p_-,p_+\in\partial X,
        \qquad
        q_-,q_+\in\partial(X+w),
\]
and pairwise disjoint open neighborhoods $P_-,P_+,V_-,V_+$ of these four points such that:
\begin{enumerate}[label=\textup{(\roman*)},leftmargin=2.2em]
\item $\mathcal N_{\varepsilon_0}(X)\subset Q$;
\item whenever $|w'-w|<\eta_0$,
\[
        \overline{P_-}\cup\overline{P_+}
        \quad\text{is disjoint from }Q+w',
        \qquad
        \overline{V_-}\cup\overline{V_+}
        \quad\text{is disjoint from }Q;
\]
\item for every $|w'-w|<\eta_0$, every $p'_-\in P_-$, $p'_+\in P_+$, $q'_-\in V_-$ and $q'_+\in V_+$, each path in $Q$ from $p'_-$ to $p'_+$ meets each path in $Q+w'$ from $q'_-$ to $q'_+$.
\end{enumerate}
\end{proposition}

\begin{proof}
See \cref{app:linked-translates}.
\end{proof}

\begin{definition}[Normalized finite trap]\label{def:finite-trap}
Let $c\in\C$ with $|c|>1$, let $w\in\C$, and put $K:=K(c,n)$.  A \emph{normalized finite trap} for $(c,w)$ consists of a compact topological disc $Q\subset\C$, a number $\varepsilon>0$, and points
\[
        p_-,p_+\in K,
        \qquad
        q_-,q_+\in K+w,
\]
such that
\begin{enumerate}[label=\textup{(\roman*)},leftmargin=2.2em]
\item $\mathcal N_{\varepsilon/2}(K)\subset Q$;
\item $p_-,p_+\notin Q+w$ and $q_-,q_+\notin Q$;
\item every polygonal path in $Q$ joining $p_-$ to $p_+$ meets every polygonal path in $Q+w$ joining $q_-$ to $q_+$;
\item $\delta_n(c)<\varepsilon$.
\end{enumerate}
The trap is \emph{adjacent} if $w=w_{u,v}(c)$ for adjacent words $u,v$ of a common positive length.  Clause \textup{(ii)} is part of the trap configuration of \cite[Definition~7.1.3]{CKW2017} and holds for the traps produced in \cref{lem:trap-extraction}; the arguments below use only \textup{(i)}, \textup{(iii)} and \textup{(iv)}.
\end{definition}

\begin{figure}[t]
\centering
\includegraphics[width=\textwidth]{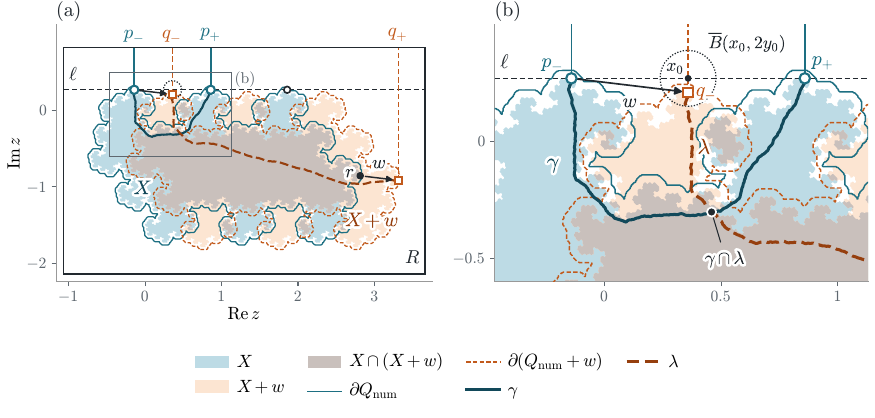}
\caption{Illustration of the stable linked translate for $X=\Fill{K(c_0,3)}$, where $c_0\approx1.0933+1.3660\,i$ is the root in the upper half-plane of $P(c)=c^4-2c^3+2c^2+2c-2$ (\cref{ex:worked}).
\textup{(a)} The upper support line $\ell$ meets $X$ in three points; $p_-,p_+$ bound its first gap and $x_0$ is their midpoint.  With $y_0=3/50$, the translation is $w=1/2-3i/50$, $q_-=p_-+w$ and $q_+=r+w$, where $r$ is the rightmost point of $X$.  The irrational direction and the disjointness of $\overline B(x_0,2y_0)$ from $X$ are proved by exact inequalities in \cref{app:worked-certificate}, so this vector satisfies the construction in \cref{app:linked-translates}.  The vertical access segments reach the upper side of $R$ in alternating order.
\textup{(b)} Enlargement of the boxed region, showing a crossing of the sampled short-hop paths.
All depth-$13$ prefixes are enumerated and evaluated numerically (paths: depth $11$).  The displayed filled hull, the filled $0.035$-neighborhood $Q_{\rm num}$ and the paths are numerical illustrations, not certified enclosures.  In particular, $Q_{\rm num}$ is not identified with the topological disc $Q$ whose existence is proved in \cref{prop:stable-linked}.}
\label{fig:finite-trap}
\end{figure}

\begin{proposition}[Adjacent finite traps force cylinder overlap]\label{prop:finite-trap}
If $c$ has an adjacent normalized finite trap with translation $w=w_{u,v}(c)$, then
\[
        K_u(c)\cap K_v(c)\ne\varnothing,
\]
equivalently $K(c,n)\cap\bigl(K(c,n)+w\bigr)\ne\varnothing$.  In particular $c\in\Mn$.
\end{proposition}

The first stage of the proof is the argument of \cite[Proposition~7.1.6]{CKW2017}; the second stage extracts the overlap of the two cylinders, which is what the perturbative criterion below uses.

\begin{proof}
Let $K:=K(c,n)$ and put $m:=|u|=|v|$.  We first show that $\delta_n(c)=0$.  Suppose instead that $\delta:=\delta_n(c)>0$, and choose $\tau>0$ so that
\[
        \delta+\tau<\varepsilon,
        \qquad
        \delta+\tau<|c|^m\delta.
\]
By \Cref{lem:short-hop}, there is a polygonal path $\gamma$ from $p_-$ to $p_+$, with vertices in $K$ and edge lengths below $\delta+\tau$.  Every edge lies in $\mathcal N_{(\delta+\tau)/2}(K)\subset Q$.  The translated version of the same lemma gives a path $\lambda\subset Q+w$ from $q_-$ to $q_+$ with vertices in $K+w$ and the same edge bound.

The crossing condition forces $\gamma\cap\lambda\ne\varnothing$.  At an intersection point, the distances to $K$ and $K+w$ are both below $(\delta+\tau)/2$, and therefore
\[
        \dist(K,K+w)<\delta+\tau.
\]
On the other hand, \cref{lem:cylinder-identities} gives
\[
        \dist(K,K+w)=|c|^m\dist\bigl(K_u(c),K_v(c)\bigr).
\]
The two cylinders lie in neighboring first-level pieces, so their distance is at least $\delta$.  Thus $\dist(K,K+w)\ge|c|^m\delta>\delta+\tau$, a contradiction.  Hence $\delta_n(c)=0$, which already gives $c\in\Mn$ by \cref{prop:neighbor}; the rest of the proof yields the overlap of the two specific cylinders.

Now let $0<\tau<\varepsilon$.  Since $\delta_n(c)=0$, \cref{lem:short-hop} gives paths $\gamma\subset\mathcal N_{\tau/2}(K)\subset Q$ from $p_-$ to $p_+$ and $\lambda\subset\mathcal N_{\tau/2}(K+w)\subset Q+w$ from $q_-$ to $q_+$, with vertices in $K$ and $K+w$ respectively and edges shorter than $\tau$.  They meet by the crossing condition, and at a common point both distances are below $\tau/2$; hence $\dist(K,K+w)<\tau$.  As $\tau$ is arbitrary and both sets are compact, $K\cap(K+w)\ne\varnothing$, and \cref{lem:cylinder-identities} gives $K_u(c)\cap K_v(c)\ne\varnothing$.
\end{proof}

\begin{lemma}[Finite trap extraction]\label{lem:trap-extraction}
Let $c_0\in\Mn$ and suppose that
\[
        X_0:=\Fill{K(c_0,n)}
\]
is nonconvex.  Then there exist a vector $w\in\C$, a neighborhood $U$ of $c_0$, a number $\eta>0$, a compact topological disc $Q$, a number $\varepsilon>0$, and four addresses
\[
        \alpha_-,\alpha_+,\beta_-,\beta_+\in\An^{\N}
\]
such that, for every $c\in U$ and every $w'\in B(w,\eta)$, the points
\[
        p_\pm(c):=\pi_c(\alpha_\pm),
        \qquad
        q_\pm(c,w'):=\pi_c(\beta_\pm)+w'
\]
together with $Q$ and $\varepsilon$ form a normalized finite trap for $(c,w')$.
\end{lemma}

\begin{proof}
Since $c_0\in\Mn$, the attractor $K(c_0,n)$ is a continuum, and its filled hull $X_0$ is a full continuum: it is compact, it is connected because each bounded complementary component of $K(c_0,n)$ has its boundary in $K(c_0,n)$, and its complement is the unbounded component of $\C\setminus K(c_0,n)$.  Apply \Cref{prop:stable-linked} to $X_0$.  By \cref{lem:filled-boundary}, the points $p_\pm$ belong to $K(c_0,n)$ and the points $q_\pm-w$ belong to $K(c_0,n)$.  Choose addresses $\alpha_\pm$ and $\beta_\pm$ for these four points.

Choose $0<\varepsilon<\varepsilon_0/2$.  The address maps are continuous in $c$, and $K(c,n)$ depends continuously on $c$ in the Hausdorff metric.  Moreover, $\delta_n(c)$ is continuous and vanishes at $c_0$ by \cref{prop:neighbor}.  After shrinking a neighborhood $U$ of $c_0$ and choosing $0<\eta<\eta_0$, we may therefore arrange that
\begin{enumerate}[label=\textup{(\alph*)},leftmargin=2.2em]
\item $p_\pm(c)\in P_\pm$ and $q_\pm(c,w')\in V_\pm$;
\item $d_{\mathrm H}(K(c,n),K(c_0,n))<\varepsilon_0/4$;
\item $\delta_n(c)<\varepsilon$
\end{enumerate}
for all $c\in U$ and $w'\in B(w,\eta)$.  Since $K(c_0,n)\subset X_0$ and $\varepsilon<\varepsilon_0/2$, condition \textup{(b)} gives
\[
        \mathcal N_{\varepsilon/2}(K(c,n))
        \subset \mathcal N_{\varepsilon_0}(X_0)
        \subset Q.
\]
The separation and crossing clauses in \cref{def:finite-trap} now follow from \cref{prop:stable-linked}; clause \textup{(iv)} is \textup{(c)}.
\end{proof}

\Needspace{9\baselineskip}
\begin{theorem}[Perturbative trap criterion]\label{thm:robust-trap}
Let $c_0\in\Mn$ and assume that $\Fill{K(c_0,n)}$ is nonconvex (so $c_0\notin\R$, a real connected attractor being an interval).  Then there exist a vector $w\in\C$, a neighborhood $U$ of $c_0$, and $\eta>0$ such that whenever $c\in U$ and adjacent words $u,v$ of a common positive length satisfy
\[
        |w_{u,v}(c)-w|<\eta,
\]
there is a neighborhood $U_c\subset U$ of $c$ with
\[
        K_u(c')\cap K_v(c')\ne\varnothing
        \qquad\text{for every }c'\in U_c;
\]
in particular $U_c\subset\Mn$ and $c\in\interior(\Mn)$.  Any vector $w$ supplied with the accompanying data of \cref{prop:stable-linked} for this filled hull may be used.
\end{theorem}

\begin{proof}
Apply \Cref{lem:trap-extraction} to $c_0$ and write $w,U,\eta$ for the resulting data.

Fix $c\in U$ and adjacent words $u,v$ satisfying the displayed inequality.  By continuity of $c'\mapsto w_{u,v}(c')$, the same strict inequality holds throughout a smaller neighborhood $U_c\subset U$ of $c$.  For every $c'\in U_c$, \cref{lem:trap-extraction} supplies an adjacent normalized finite trap for $(c',w_{u,v}(c'))$, and \cref{prop:finite-trap} gives $K_u(c')\cap K_v(c')\ne\varnothing$ and $c'\in\Mn$.  Hence $U_c\subset\Mn$ and $c\in\interior(\Mn)$.
\end{proof}

\Cref{lem:trap-extraction,thm:robust-trap} are the collinear form of the stability step in the proof of \cite[Theorem~7.2.7]{CKW2017}.  Two points differ from the binary case.  First, the neighboring-piece distance $\delta_n(c)$ replaces the distance between the two pieces of the binary attractor; that it vanishes exactly on $\Mn$, and that its vanishing propagates to chains with short hops through all $n$ pieces, is the content of \cref{prop:neighbor,lem:short-hop}.  Second, the conclusion records the overlap of the two specific cylinders $K_u,K_v$ throughout a neighborhood, rather than membership in $\interior(\Mn)$ alone; \cref{ex:worked} uses this form.

\section{Proofs of the regular-closedness theorems}\label{sec:main-proof}

\begin{proof}[Proof of \cref{thm:main}]
Fix $c_0\in\Mn\setminus\R$.

If $1<|c_0|<\sqrt n$, then $c_0\in\interior(\Mn)$ by \cref{lem:annulus}.  If $|c_0|=\sqrt n$, the radial parameters
\[
        c_r:=\frac{r}{|c_0|}c_0,
        \qquad 1<r<\sqrt n,
\]
belong to $\interior(\Mn)$ and converge to $c_0$ as $r\uparrow\sqrt n$.  Thus $c_0\in\overline{\interior(\Mn)}$ in both cases.

Assume $|c_0|>\sqrt n$.  By \Cref{cor:convex-bound}, the filled hull $\Fill{K(c_0,n)}$ is nonconvex.  Let $W$ be an arbitrary neighborhood of $c_0$.  Apply \Cref{thm:robust-trap} and shrink its neighborhood $U$ so that $U\subset W$.

By \Cref{prop:neighbor}, choose addresses $\alpha,\beta\in\An^{\N}$ with
\[
        |\alpha_0-\beta_0|=1,
        \qquad
        \pi_{c_0}(\alpha)=\pi_{c_0}(\beta).
\]
Apply \Cref{lem:rouche} with the target vector $w$ furnished by \cref{thm:robust-trap} and the neighborhood $U$.  We obtain a parameter $c\in U$ and adjacent prefixes $u,v$ such that $w_{u,v}(c)=w$.  The perturbative trap criterion yields $c\in\interior(\Mn)$.  Since $W$ was arbitrary, $c_0\in\overline{\interior(\Mn)}$.

This proves
\[
        \Mn\setminus\R\subset\overline{\interior(\Mn)}.
\]
The reverse inclusion in \eqref{eq:main}, after removing the real axis, follows from the closedness of $\Mn$ in \cref{cor:power-series}.  Finally, non-real points are dense in the open set $\interior(\Mn)$, so $\overline{\interior(\Mn)\setminus\R}=\overline{\interior(\Mn)}$, and \eqref{eq:main} is the same as regular closedness of $\Mn\setminus\R$ in $\{c\in\C\setminus\R:|c|>1\}$.
\end{proof}

The next consequence relates the geometric statement to polynomial roots outside the closed unit disc.  In this corollary, closures and interiors remain relative to $\{z\in\C:|z|>1\}$.

\begin{corollary}[Polynomial roots and the regular core]\label{cor:regular-core}
Let $\mathcal R_n$ be the set of roots of modulus greater than $1$ of monic integer polynomials of height at most $n-1$, that is, of the polynomials
\[
        P(z)=z^m+a_1z^{m-1}+\cdots+a_m,
        \qquad m\ge1,\quad a_j\in\mathbb Z,\quad |a_j|\le n-1.
\]
Then
\begin{equation}\label{eq:regular-core}
\begin{gathered}
        \Mn=\overline{\mathcal R_n},\qquad
        \Core=\overline{\Mn\setminus\R}
              =\overline{\mathcal R_n\setminus\R}
              =\overline{\mathcal R_n\cap\interior(\Mn)},\\
        \interior(\Core)=\interior(\Mn).
\end{gathered}
\end{equation}
Moreover,
\begin{equation}\label{eq:decomposition}
        \Mn=\Core\cup(\partial\Mn\cap\R),
        \qquad
        \Mn\setminus\Core\subset\partial\Mn\cap\R,
\end{equation}
and the non-real boundary of the locus is the boundary of its interior:
\begin{equation}\label{eq:boundary}
        \partial\Mn\setminus\R
        =\partial\Core\setminus\R
        =\partial\,\interior(\Mn)\setminus\R.
\end{equation}
\end{corollary}

The identity $\Mn=\overline{\mathcal R_n}$ is \cite[Proposition~2.3]{EJ2026FiniteCapture}; the argument is included because the non-real refinements in \eqref{eq:regular-core} depend on the choice of the Rouch\'e circles made below.

\begin{proof}
Padding the coefficients with zeros in \eqref{eq:power-series} gives
$\mathcal R_n\subset\Mn$.  Conversely, for $c_0\in\Mn$, choose an
admissible series $F(c)=1-\sum_{k\ge0}d_kc^{-(k+1)}$ vanishing at $c_0$.
Its truncations converge locally uniformly on $\{z\in\C:|z|>1\}$, and $F(c)\to1$
as $|c|\to\infty$, so $F$ is not identically zero.
Rouch\'e's theorem on arbitrarily small isolating circles
about $c_0$ gives nearby zeros of these truncations.  Multiplication by
the corresponding power of $c$ makes each truncation a monic polynomial
of the stated form.  When $c_0$ is non-real, the circles may be chosen
disjoint from $\R$.  Thus
$\overline{\mathcal R_n}=\Mn$ and
$\overline{\mathcal R_n\setminus\R}=\overline{\Mn\setminus\R}$.

Non-real points are dense in the open set $\interior(\Mn)$, so
$\Core\subset\overline{\Mn\setminus\R}$.  The reverse inclusion follows
from \cref{thm:main}.  If $c\in\interior(\Mn)$ and $V\subset\Mn$ is an open
neighborhood of $c$, the roots approximating $c$ eventually lie in $V$ and
hence in $\interior(\Mn)$; thus
$\interior(\Mn)\subset\overline{\mathcal R_n\cap\interior(\Mn)}\subset\Core$,
and taking closures gives the last identity in the first line of
\eqref{eq:regular-core}.  Next,
$\interior(\Mn)\subset\Core\subset\Mn$ gives equality of interiors.

Now $\Mn$ is closed in $\{z\in\C:|z|>1\}$ by \cref{cor:power-series}, so
$\partial\Mn=\Mn\setminus\interior(\Mn)$ and $\Core\cup(\partial\Mn\cap\R)\subset\Mn$.
Conversely, a point of $\Mn\setminus\Core$ is real by \cref{thm:main} and is
not interior because $\interior(\Mn)\subset\Core$; hence it lies in
$\partial\Mn\cap\R$.  This proves \eqref{eq:decomposition}.

Finally, $\Core$ is closed, so by \eqref{eq:regular-core}
$\partial\Core=\Core\setminus\interior(\Core)=\Core\setminus\interior(\Mn)
=\partial\,\interior(\Mn)$, while $\partial\Mn=\Mn\setminus\interior(\Mn)$.
Since $\Mn\setminus\R=\Core\setminus\R$ by \cref{thm:main}, removing the real
axis gives \eqref{eq:boundary}.
\end{proof}

We now pass from $\Mn\subset\{c\in\C:|c|>1\}$ to the full planar root closure.  The next proposition is the identification of \cite[Proposition~1.1]{EJ2026FiniteCapture}, restricted to the plane.  Its proof is included to make the deduction of \cref{thm:root-closure} independent of that reference.  All operations on $\mathcal Z_h$ are taken in $\C$.

\begin{proposition}[Full-plane root-closure identification]\label{prop:planar-identification}
Let $h\ge1$, put $n=h+1$, and let $\iota(z)=1/z$ for $z\ne0$.  Then
\begin{equation}\label{eq:planar-root-decomposition}
        \mathcal Z_h
        =\Mn\cup\iota(\Mn)\cup\{z\in\C:|z|=1\}\cup\{0\}.
\end{equation}
The same planar closure is obtained by restricting to monic integer polynomials of height at most $h$.
\end{proposition}

\begin{proof}
Let $P(z)=\sum_{j=0}^d a_jz^j$ have height at most $h$ and let $P(c)=0$ with $|c|>1$.  Changing the sign of $P$ if necessary, assume $1\le a_d\le h$.  Write
\[
        c^{-d}P(c)=\sum_{j=0}^d (u_j-v_j)c^{-j},
        \qquad u_j,v_j\in\{0,\ldots,h\},
\]
by taking $u_j=\max\{a_{d-j},0\}$ and $v_j=\max\{-a_{d-j},0\}$.  Padding both words with zeros gives two addresses of the same point, with distinct initial digits $u_0=a_d$ and $v_0=0$.  Thus two distinct first-level pieces intersect.  The final assertion of \cref{prop:neighbor} forces $c\in\Mn$.  Closedness of $\Mn$ and \cref{cor:regular-core} now give
\[
        \mathcal Z_h\cap\{|z|>1\}=\Mn
        =\overline{\mathcal R_n},
\]
where the last closure is relative to $\{|z|>1\}$.

The set of nonzero polynomial roots under consideration is invariant under inversion: after removing any factor $z^k$ from $P$, reverse the remaining coefficients.  Reversal preserves height and sends a nonzero root to its reciprocal.  Since inversion is a homeomorphism of $\C\setminus\{0\}$, the same invariance holds for $\mathcal Z_h\setminus\{0\}$.  Roots of $z^m-1$ are dense on the unit circle, and $0$ is a root of $z$.  These observations prove \eqref{eq:planar-root-decomposition}.

It remains to verify the monic assertion inside the unit disc.  If $a$ with $|a|<1$ is a root of a nonzero polynomial $P$ of height at most $h$, then, for $N>\deg P$, the polynomial $z^N+P(z)$ is monic and has height at most $h$.  On a sufficiently small isolating circle about $a$ contained in $\{|z|<1\}$, $z^N+P(z)\to P(z)$ uniformly.  Rouch\'e's theorem gives roots of these monic polynomials converging to $a$.  Outside the closed unit disc the monic assertion follows from $\Mn=\overline{\mathcal R_n}$, and the roots used on the unit circle and at $0$ are already monic.  Taking planar closures proves the claim.
\end{proof}

\begin{proof}[Proof of \cref{thm:root-closure}]
Put $n=h+1$ and let $\iota(z)=1/z$ on $\C\setminus\{0\}$.  We use the decomposition \eqref{eq:planar-root-decomposition} proved in \cref{prop:planar-identification}.

Let $z\in\mathcal Z_h\setminus\R$.  If $|z|>1$, then $z\in\Mn$ by \eqref{eq:planar-root-decomposition}.  By \cref{thm:main}, it is a limit of points of $\interior(\Mn)$.  Each such point lies in $\interior_{\C}(\mathcal Z_h)$: the relative interior $\interior(\Mn)$ is open in the open set $\{c\in\C:|c|>1\}$, hence open in $\C$, and $\Mn\subset\mathcal Z_h$.

If $0<|z|<1$, then $z\in\iota(\Mn)$ by \eqref{eq:planar-root-decomposition}, so $1/z\in\Mn\setminus\R$.  Choose $c_j\in\interior(\Mn)$ with $c_j\to1/z$ by \cref{thm:main}.  Inversion is a homeomorphism of $\C\setminus\{0\}$, so $1/c_j$ lies in the planar interior of $\iota(\Mn)\subset\mathcal Z_h$, and $1/c_j\to z$.

Finally, \cref{lem:annulus} and \eqref{eq:planar-root-decomposition} imply
\begin{equation}\label{eq:unit-circle-annulus}
        \{z\in\C:n^{-1/2}<|z|<n^{1/2}\}
        \subset\interior_{\C}(\mathcal Z_h).
\end{equation}
Indeed, its part outside the unit circle lies in $\Mn$, its part inside lies in $\iota(\Mn)$, and the unit circle itself is included in \eqref{eq:planar-root-decomposition}; the whole annulus is open.  Thus points with $|z|=1$ are already interior.  The point $0$ is real and does not occur in the assertion.  These cases prove the forward inclusion in \eqref{eq:global-regular-closedness}; the reverse follows because $\mathcal Z_h$ is closed in $\C$.  Since non-real points are dense in the open set $\interior_{\C}(\mathcal Z_h)$, the right-hand side of \eqref{eq:global-regular-closedness} is the closure in $\C\setminus\R$ of $\interior_{\C}(\mathcal Z_h)\setminus\R$, so \eqref{eq:global-regular-closedness} is the same as regular closedness of $\mathcal Z_h\setminus\R$ in $\C\setminus\R$.
\end{proof}

\begin{example}[A worked instance of the proof]\label{ex:worked}
Let $n=3$ and let $c_0$ be the root in the upper half-plane of
\[
        P(c)=c^4-2c^3+2c^2+2c-2,
        \qquad c_0\approx1.0933+1.3660\,i.
\]
The rational isolation in \cref{app:worked-certificate} proves that this root is simple, that $|c_0|>\sqrt3$, and that $|c_0+1|>\sqrt6$.  Thus $c_0$ lies outside the lens $\mathcal X_3$ of \cref{sec:structure}.  Dividing $P(c_0)=0$ by $c_0^4$ gives \eqref{eq:power-series} with $(d_0,d_1,d_2,d_3)=(2,-2,-2,2)$ and $d_k=0$ for $k\ge4$, so $c_0\in\mathcal M_3\setminus\R$ by \cref{cor:power-series}.  The same identity gives equal address values for
\[
        \alpha=0\,2\,0\,0\,2\,0^\infty,
        \qquad
        \beta=1\,0\,2\,2\,0^\infty.
\]
These addresses are adjacent, and in the notation of \cref{lem:rouche},
$T(c)=-P(c)c^{-4}$ and $T_m=T$ for $m\ge5$.

Fix
\[
        w=\frac12-\frac{3}{50}i.
\]
\Cref{app:worked-certificate} proves that this particular vector satisfies the supporting-line construction of \cref{app:linked-translates} for $X_0=\Fill{K(c_0,3)}$.  Applying the extraction and stability arguments with this vector gives $U$ and $\eta$ as in \cref{thm:robust-trap}.  Figure~\ref{fig:finite-trap} illustrates the configuration; no numerical enclosure is needed for this conclusion.

For $m\ge5$, the realization equation is
\begin{equation}\label{eq:worked-realization}
        c^{m-4}P(c)+w=0.
\end{equation}
On a fixed sufficiently small disc about $c_0$ contained in $\{|c|>1\}$, Rouch\'e's theorem applied to $P(c)+wc^{4-m}$ gives a unique root $c_m$ for every sufficiently large $m$, and $c_m\to c_0$.  It satisfies
$w_{\alpha|_m,\beta|_m}(c_m)=w$ exactly.  Moreover,
\begin{equation}\label{eq:worked-asymptotic}
        c_m-c_0=-\frac{w\,c_0^{\,4-m}}{P'(c_0)}
                  +o\!\left(|c_0|^{-m}\right),
        \qquad
        \frac{|c_{m+1}-c_0|}{|c_m-c_0|}\longrightarrow |c_0|^{-1}.
\end{equation}
Indeed, for some fixed $r>1$, \eqref{eq:worked-realization} and the simple zero of $P$ give $c_m-c_0=O(r^{-m})$.  Hence $m(c_m-c_0)\to0$ and $c_m^{4-m}/c_0^{4-m}\to1$; Taylor expansion of $P$ gives \eqref{eq:worked-asymptotic}.  Nearby solutions at depths $8$ and $12$ are numerically $1.0899+1.3706\,i$ and $1.0935+1.3655\,i$, respectively; these values are illustrative and are not finite-depth interior certificates.

For all sufficiently large $m$, $c_m\in U$, and the fixed-vector form of \cref{thm:robust-trap} gives a neighborhood of $c_m$ on which $K_{\alpha|_m}$ and $K_{\beta|_m}$ overlap.  In particular, $c_m\in\interior(\mathcal M_3)$.  The proof thus supplies interior parameters accumulating at $c_0$; it does not establish that $c_0$ itself is interior.
\end{example}

\section{Internal structure of the locus}\label{sec:structure}

\Cref{thm:main} and \cref{cor:regular-core} express $\Mn$ through the regular core $\Core$ and the real boundary trace $\partial\Mn\cap\R$; we call $\Mn\cap\R$ the \emph{real trace} and $\Mn\cap i\R$ the \emph{imaginary trace}.  This section records what the preceding results, together with the capture levels that Espigul\'e and Juher construct inside a lens \cite{EJ2026FiniteCapture}, say about these two pieces.  Closures, interiors and boundaries remain relative to $\{z\in\C:|z|>1\}$; this is also the ambient convention of \cite{EJ2026FiniteCapture}.

\subsection{The real and imaginary traces}

\begin{proposition}[Real and imaginary traces of the locus]\label{prop:real-trace}
Let $n\ge2$.
\begin{enumerate}[label=\textup{(\roman*)},leftmargin=2.2em]
\item $\Mn\cap\R=\{\sigma\in\R:1<|\sigma|\le n\}$.
\item $\{\sigma\in\R:1<|\sigma|<\sqrt n\}\subset\interior(\Mn)$, and $\pm\sqrt n\in\Core$.
\item Every real $\sigma$ with $1+\sqrt{n-1}<|\sigma|\le n$ belongs to $\partial\Mn\setminus\Core$.
\item $\Mn\cap i\R=\{iy:y\in\R,\ 1<|y|\le\sqrt n\}$.  All these points belong to $\Core$, and $\pm i\sqrt n\in\partial\Mn$.
\end{enumerate}
Consequently,
\begin{equation}\label{eq:real-sandwich}
        \{\sigma\in\R:1+\sqrt{n-1}<|\sigma|\le n\}
        \subset\Mn\setminus\Core
        \subset\{\sigma\in\R:\sqrt n<|\sigma|\le n\}.
\end{equation}
\end{proposition}

\begin{proof}
\textup{(i)} Let $\sigma\in\R$ with $|\sigma|>1$.  The attractor $K(\sigma,n)$ is a compact subset of $\R$ containing the two distinct points $\pi_\sigma(0^\infty)=0$ and $\pi_\sigma((n-1)^\infty)=(n-1)\sigma/(\sigma-1)$; hence it is connected if and only if it is a nondegenerate interval.  By \eqref{eq:centered}, $K(\sigma,n)$ is a translate of $\Ktil(\sigma,n)$, and \cref{lem:real-criterion} shows that $\Ktil(\sigma,n)$ is an interval when $|\sigma|\le n$ and contains no nondegenerate interval when $|\sigma|>n$.

\textup{(ii)} The first assertion is \cref{lem:annulus}.  For the second, the real parameters $\pm r$ with $1<r<\sqrt n$ are interior by the first assertion and converge to $\pm\sqrt n$ as $r\uparrow\sqrt n$.

\textup{(iii)} By \cite[Proposition~A.2]{EJ2026FiniteCapture},
\begin{equation}\label{eq:nonreal-bound}
        \Mn\setminus\R\subset\{c\in\C:|c|<1+\sqrt{n-1}\};
\end{equation}
the closed-disc version is \cite[Proposition~4]{EJS2024}, derived there from a radial bound of Beaucoup, Borwein, Boyd and Pinner \cite{BBBP1998}.  Let $\sigma\in\R$ with $1+\sqrt{n-1}<|\sigma|\le n$.  Then $\sigma\in\Mn$ by \textup{(i)}, and
\[
        V:=\{c\in\C:|c|>1+\sqrt{n-1}\}
\]
is an open neighborhood of $\sigma$ with $V\cap\Mn\subset\R$ by \eqref{eq:nonreal-bound}.  If $V$ met $\interior(\Mn)$, then $V\cap\interior(\Mn)$ would be a nonempty open subset of $\Mn$, and every nonempty open subset of $\C$ contains non-real points; this contradicts $V\cap\Mn\subset\R$.  Hence $V\cap\interior(\Mn)=\varnothing$, so $\sigma\notin\overline{\interior(\Mn)}=\Core$, while $\sigma\in\Mn\setminus\interior(\Mn)=\partial\Mn$.

\textup{(iv)} Let $c=iy$ with $y\in\R$ and $|y|>1$, so that $c^2=-y^2$.  Grouping the series in \eqref{eq:centered} according to the parity of the index, as in the proof of \cref{thm:convex-locus}, gives
\[
        \Ktil(iy,n)=A+(iy)^{-1}A=\{a-ia'/y:a,a'\in A\},
        \qquad
        A:=\Ktil(-y^2,n)\subset\R.
\]
Thus $\Ktil(iy,n)$ is the image of $A\times A$ under a real-linear isomorphism of $\R^2$ onto $\C$, and it is connected if and only if $A$ is.  The compact set $A$ is not a single point, its convex hull being $[-M,M]$ with $M>0$ as in the proof of \cref{lem:real-criterion}; hence $A$ is connected if and only if it is a nondegenerate interval, which by \cref{lem:real-criterion} happens exactly when $y^2\le n$.  This proves the description of $\Mn\cap i\R$.  Its points with $|y|<\sqrt n$ are interior by \cref{lem:annulus}, so all of its points lie in $\Core$, and $\pm i\sqrt n$ are limits of the points $iy\notin\Mn$ with $|y|\downarrow\sqrt n$, hence belong to $\partial\Mn$.

The first inclusion in \eqref{eq:real-sandwich} is \textup{(iii)}.  For the second, $\Mn\setminus\Core\subset\R$ by \cref{thm:main}, real parameters with $|\sigma|<\sqrt n$ are interior and $\pm\sqrt n\in\Core$ by \textup{(ii)}, and $|\sigma|\le n$ by \textup{(i)}.
\end{proof}

\begin{remark}[Real whiskers]\label{rem:whiskers}
For $n=2$ the set in \textup{(iii)} is empty, since $1+\sqrt{n-1}=n$.  In that case more is known: in the parameter $z=1/c$, Calegari, Koch and Walker showed that there are $t_*>1/2$ and an open set $U\subset\C$ with $\{1/c:c\in\mathcal M_2\}\cap U=[1/2,t_*)$ \cite[Lemma~10.1.1]{CKW2017}.  The argument in \textup{(iii)} then shows that the real parameters $c$ with $1/t_*<c\le2$ lie in $\partial\mathcal M_2\setminus\mathcal K_2$; see \cite[Section~10]{CKW2017} for the finer structure near the real axis.  The full real trace of the regular core in the range $\sqrt n<|\sigma|\le1+\sqrt{n-1}$ is not determined by these results.  It is not empty in general, as the following example shows.
\end{remark}

\begin{example}[A real core point beyond the interior annulus]\label{ex:real-core}
For $n=5$, set $\sigma=1+\sqrt2$.  Then
\[
        \sqrt5<\sigma<3=1+\sqrt{5-1},
        \qquad \sigma\in\mathcal K_5.
\]
Indeed, let
\[
        F(z)=(z^2-2z-1)^2=z^4-4z^3+2z^2+4z+1,
        \qquad P_m(z)=z^mF(z)+1\quad(m\ge1).
\]
Each $P_m$ is monic of height at most $4$.  On a small closed disc about $\sigma$ contained in $\{|z|>1,\ \Re z>0\}$ and containing no other zero of $F$, the functions $F(z)+z^{-m}$ converge uniformly to $F(z)$.  Rouch\'e's theorem gives two zeros, counted with multiplicity, of $P_m$ in this disc for all sufficiently large $m$.  Shrinking the disc shows that these zeros converge to $\sigma$.  They are non-real because
\[
        P_m(x)=x^m(x^2-2x-1)^2+1>0\qquad(x>0).
\]
Thus $\sigma\in\overline{\mathcal R_5\setminus\R}=\mathcal K_5$ by \cref{cor:regular-core}.  Substitution $z\mapsto-z$, followed by a sign change when necessary to keep the polynomial monic, also gives $-\sigma\in\mathcal K_5$.
\end{example}

\subsection{Capture levels and the regular core}

For $n\ge2$ let
\begin{equation}\label{eq:lens}
        \Xn:=\{c\in\C:|c|>1,\ |c-1|<\sqrt{2n}\ \text{and}\ |c+1|<\sqrt{2n}\}
\end{equation}
be the open lens of \cite[(1.4)]{EJ2026FiniteCapture}; we call $\Xn\setminus\R$ the \emph{non-real lens}.  Espigul\'e and Juher use the symmetric alphabet of \cref{lem:normalization}, under which the difference attractor is $E(c,N)=2D(c,n)$, where $N:=2n-1$, and the marked point is $2c$; the lens and the sets defined below lie in parameter space and do not depend on this choice.  Call a nonempty bounded open set \emph{self-covering} if it is contained in the union of its images under the $N$ contractions $z\mapsto t+z/c$, $t\in A_N^{\mathrm{sym}}$, of the difference system.  For non-real $c$, a self-covering parallelogram with sides parallel to $1$ and $\overline c$ exists if and only if $c\in\Xn$ \cite[Theorem~A]{EJ2026FiniteCapture}.  For $c\in\Xn\setminus\R$, the \emph{canonical parallelogram} $\mathcal C(c,N)$ is such a set, the \emph{canonical trap}, and it lies in the interior of $E(c,N)$ \cite[Definition~3.3 and Corollary~3.4]{EJ2026FiniteCapture}.  A word $u$ of digits \emph{captures} $c$ if the inverse branch $g_u$ of the difference system indexed by $u$ sends $2c$ into $\mathcal C(c,N)$.  The \emph{capture level} $\Theta_k(n)\subset\Xn\setminus\R$ consists of the parameters captured by a word of length at most $k$, the \emph{first capture depth} $\kappa_n(c)$ is the least such $k$, and
\[
        \Thn:=\bigcup_{k\ge0}\Theta_k(n)
\]
is the union of all capture levels \cite[Section~1.2]{EJ2026FiniteCapture}.  As explained in \cref{sec:introduction}, the canonical trap is not a trap in the sense of \cref{sec:traps}.  We use the following properties.

\Needspace{9\baselineskip}
\begin{theorem}[Capture levels; Espigul\'e--Juher]\label{thm:finite-capture}
Let $n\ge2$.
\begin{enumerate}[label=\textup{(\alph*)},leftmargin=2.2em]
\item The capture levels increase with $k$; each $\Theta_k(n)$ is open and semialgebraic, and $\Theta_k(n)\subset\interior(\Mn)$ \cite[Theorem~B and Proposition~4.1]{EJ2026FiniteCapture}.  Hence $\Thn$ is open and $\Thn\subset\interior(\Mn)$.
\item A non-real root $c\in\Xn$ of a polynomial $z^m+a_1z^{m-1}+\cdots+a_m$ as in \cref{cor:regular-core} lies in $\Theta_m(n)$: the digits $-2a_1,\ldots,-2a_m$ form a word of length $m$ whose inverse branch sends $2c$ to $0\in\mathcal C(c,N)$ \cite[Proposition~2.4]{EJ2026FiniteCapture}.  In particular $\Rn\cap(\Xn\setminus\R)\subset\Thn$.
\item $\overline{\Thn}\cap(\Xn\setminus\R)=\Mn\cap(\Xn\setminus\R)$ \cite[Theorem~D]{EJ2026FiniteCapture}.
\item If $n\ge20$, then $\Mn\setminus\R\subset\Xn$; for each $2\le n\le19$ there is a non-real point of $\Mn$ outside $\Xn$ \cite[Theorem~5.3]{EJ2026FiniteCapture}.
\item If $n\ge3$, the parameter $\gamma_n:=(1+i\sqrt{4n-1})/2$ belongs to $\interior(\Mn)\cap(\Xn\setminus\R)$ but not to $\Thn$, while every other parameter in some neighborhood of $\gamma_n$ belongs to $\Thn$ \cite[Theorem~5.2]{EJ2026FiniteCapture}.
\end{enumerate}
\end{theorem}

\begin{corollary}[The regular core and the capture levels]\label{cor:finite-capture}
For every $n\ge2$,
\begin{equation}\label{eq:core-capture}
        \Core\cap\Xn=\overline{\Thn}\cap\Xn,
        \qquad
        \Rn\cap(\Xn\setminus\R)\subset\Thn\subset\interior(\Mn).
\end{equation}
If $n\ge20$, then
\begin{equation}\label{eq:core-capture-global}
        \Core=\overline{\Thn},
        \qquad
        \Mn=\overline{\Thn}\cup(\partial\Mn\cap\R).
\end{equation}
For $2\le n\le19$, the off-lens set $\Mn\setminus(\R\cup\Xn)$ is nonempty and is contained in $\Core$.
\end{corollary}

\begin{proof}
The inclusions in \eqref{eq:core-capture} are \cref{thm:finite-capture}\textup{(b)} and \textup{(a)}, and \textup{(a)} gives $\overline{\Thn}\subset\Core$.  Conversely, let $c\in\Core\cap\Xn$.  By \cref{cor:regular-core}, $c$ is a limit of non-real roots $c_j\in\Rn$; since $\Xn$ is open, $c_j\in\Xn\setminus\R$ for large $j$, so $c_j\in\Thn$ by \cref{thm:finite-capture}\textup{(b)} and $c\in\overline{\Thn}$.  This proves the identity in \eqref{eq:core-capture}; off the real axis it also follows from \cref{thm:main} and \cref{thm:finite-capture}\textup{(c)}.

Let $n\ge20$.  By \cref{thm:finite-capture}\textup{(d)} and \textup{(c)}, $\Mn\setminus\R=\Mn\cap(\Xn\setminus\R)\subset\overline{\Thn}$, so $\Core=\overline{\Mn\setminus\R}\subset\overline{\Thn}$ by \cref{cor:regular-core}.  Conversely, $\Thn\subset\interior(\Mn)$ gives $\overline{\Thn}\subset\Core$.  The second identity in \eqref{eq:core-capture-global} now follows from \eqref{eq:decomposition}.  This recovers both identities of \cite[Corollary~5.4]{EJ2026FiniteCapture}.

For $2\le n\le19$, nonemptiness is \cref{thm:finite-capture}\textup{(d)} and containment is \cref{thm:main}.
\end{proof}

\begin{remark}[Capture depth, interior deficits, and holes]\label{rem:structure}
Boundaries of capture levels are transient: every limit in $\Xn\setminus\R$ of parameters in $\Theta_k(n)$ lies in $\Theta_{k+2}(n)$, and the delay two is optimal already for $k=0$ \cite[Theorem~C and Example~4.4]{EJ2026FiniteCapture}.  By \cref{cor:finite-capture} the capture levels exhaust the regular core inside the lens up to closure, but they do not exhaust the interior.  For $n\ge3$, the parameter $\gamma_n$ of \cref{thm:finite-capture}\textup{(e)}, a root of $c^2-c+n$ on the circle $|c|=\sqrt n$ bounding the annulus of \cref{lem:annulus}, is interior and belongs to no capture level, while the first capture depth of nearby parameters grows like $\log(1/|c-\gamma_n|)/\log\sqrt n$ up to a bounded error \cite[Theorem~5.2]{EJ2026FiniteCapture}.  Thus $\Thn$ is a proper subset of $\interior(\Mn)\cap(\Xn\setminus\R)$ for $n\ge3$, and \eqref{eq:core-capture} recovers the interior only up to closure.

Outside the lens, which matters precisely for $2\le n\le19$ by \cref{thm:finite-capture}\textup{(d)}, \cref{thm:main} places $\Mn\setminus\R$ in the core in two ways.  Where $|c|\le\sqrt n$, which off the lens happens only for $n\le5$, \cref{lem:annulus} suffices.  Elsewhere the mechanism is the perturbative one of \cref{sec:traps}: the interior parameter produced by \cref{lem:rouche} is a perturbation of the given one, and by \cref{cor:regular-core} it is in turn a limit of interior polynomial roots.  Inside the lens no perturbation is needed, since every non-real root there is interior by \cref{thm:finite-capture}\textup{(a)} and \textup{(b)}.

Finally, $\Mn$ has infinitely many holes for every $n\ge2$ \cite[Theorem~1.2]{Espigule2026Holes}; for $n=2$ this is the theorem of \cite{CKW2017}, and for $n\ge3$ the holes are obtained from finite-capture loops.
\end{remark}

\begin{remark}[Finite certificates off the lens]\label{rem:off-lens}
\Cref{ex:worked} exhibits the asymmetry between the two regimes.  The parameter $c_0\in\mathcal M_3\setminus(\R\cup\mathcal X_3)$ is an explicit polynomial root, yet the proof places it in $\mathcal K_3$ only through the interior parameters $c_m\to c_0$; inside the lens, by contrast, a non-real root of degree $m$ lies in $\Theta_m(n)$ and is interior by \cref{thm:finite-capture}\textup{(a)} and \textup{(b)}.  Whether every non-real root in $\Mn\setminus\Xn$, or at least $c_0$, is an interior point remains open here.

The canonical parallelogram mechanism is unavailable off the lens: no self-covering parallelogram with sides parallel to $1$ and $\overline c$ exists there \cite[Theorem~A]{EJ2026FiniteCapture}.  Whether another finite-entry construction can describe the off-lens locus is a separate question; see \cite[Section~7.4]{EJ2026FiniteCapture}.  The present proof uses stable topological crossings instead and does not require such a construction.
\end{remark}

\appendix

\section{Stable linked translates of nonconvex full continua}\label{app:linked-translates}

This appendix proves \cref{prop:stable-linked}.  The supporting-line construction is the planar core of the trap argument of Calegari--Koch--Walker \cite[Lemma~7.2.2]{CKW2017}; the details below include the neighborhood and perturbation statements needed in \cref{lem:trap-extraction}.

\begin{lemma}[Jordan-disc neighborhoods]\label{lem:jordan-neighborhood}
Let $X\subset\C$ be a full continuum.  For every $d>0$ there are a compact topological disc $Q$ and a number $e>0$ such that
\[
        \mathcal N_e(X)\subset Q\subset\mathcal N_d(X).
\]
\end{lemma}

\begin{proof}
If $X$ is a singleton, take a closed disc centered at its point with radius strictly between $0$ and $d$, and then take a smaller positive value of $e$.  Suppose henceforth that $X$ is nondegenerate.  The domain $\widehat{\C}\setminus X$ is simply connected and contains $\infty$.  Let
\[
        \phi:\{z\in\widehat{\C}:|z|>1\}\longrightarrow\widehat{\C}\setminus X
\]
be a conformal isomorphism with $\phi(\infty)=\infty$.  For $s>1$, the curve $\phi(\{|z|=s\})$ is a Jordan curve.  Let $Q_s$ be the closure of its bounded complementary component; it is a compact topological disc by the Jordan--Schoenflies theorem.  The set $\phi(\{|z|>s\})$ is open, connected, unbounded and has frontier $\phi(\{|z|=s\})$, so it is the unbounded complementary component of that curve, and $Q_s=X\cup\phi(\{1<|z|\le s\})$.  Hence $X\subset\interior(Q_s)$, the sets $Q_s$ decrease as $s\downarrow1$, and
\[
        \bigcap_{s>1}Q_s=X.
\]
Indeed, a point of $\widehat{\C}\setminus X$ has the form $\phi(\zeta)$ with $|\zeta|>1$ and lies outside $Q_s$ whenever $1<s<|\zeta|$.

It follows by compactness that $Q_s\subset\mathcal N_d(X)$ for all $s>1$ sufficiently close to $1$.  For such an $s$, the compact set $X$ lies in the interior of $Q_s$, so $\dist(X,\C\setminus Q_s)>0$.  Taking $Q=Q_s$ and any smaller positive number for $e$ proves the claim.
\end{proof}

\begin{lemma}[Alternating endpoints]\label{lem:alternating-endpoints}
Let $R$ be a compact topological disc and let $a_1,b_1,a_2,b_2$ be distinct points of $\partial R$ in this cyclic order.  Every path in $R$ joining $a_1$ to $a_2$ meets every path in $R$ joining $b_1$ to $b_2$.
\end{lemma}

\begin{proof}
Choose a homeomorphism from $R$ to a rectangle taking $a_1,a_2$ to the two vertical sides and $b_1,b_2$ to the two horizontal sides; an order-preserving homeomorphism of the boundary extends over the disc.  The assertion is then \cref{lem:rectangle-crossing}, applied to the images of the two paths.
\end{proof}

\begin{proof}[Proof of \cref{prop:stable-linked}]
We first find a gap in $X$ on a supporting line of $\conv X$.  If $\partial\conv X\subset X$, then every point of $\conv X\setminus X$ is separated from infinity by $\partial\conv X$ and hence belongs to a bounded component of $\C\setminus X$.  Since $X$ is full, this would give $\conv X=X$, contrary to nonconvexity.  Thus choose $y\in\partial\conv X\setminus X$ and a supporting line $\ell$ of $\conv X$ through $y$.  The exposed face $\conv X\cap\ell$ is a segment whose endpoints are extreme points of $\conv X$ and therefore belong to $X$: indeed, an extreme point represented as a finite convex combination of points of $X$ must coincide with every point carrying positive weight.  Since $y\notin X$, the point $y$ is not an endpoint, so the face is nondegenerate.  The component of $(\conv X\cap\ell)\setminus X$ containing $y$ is consequently an open interval whose endpoints belong to $X$.

After an orientation-preserving Euclidean motion, we may assume that
\[
        X\subset\{z:\Im z\le0\},
        \qquad
        p_-=(a,0),\quad p_+=(b,0),\quad a<b,
\]
and that the open segment between $p_-$ and $p_+$ is disjoint from $X$.  Choose $x_0\in(a,b)$ and $y_0>0$ so small that the closed disc $\overline B(x_0,2y_0)$ is disjoint from $X$, and put
\[
        v:=x_0-iy_0,
        \qquad
        w:=v-p_-.
\]
Let $M:=\max_{z\in X}\Re z$, and choose $r\in X$ with $\Re r=M$ and maximal imaginary part among the points of $X$ on that vertical line.  Set
\[
        q_-:=v=p_-+w,
        \qquad
        q_+:=r+w.
\]
The four real parts satisfy
\begin{equation}\label{eq:access-order}
        \Re p_-<\Re q_-<\Re p_+<\Re q_+.
\end{equation}
Moreover $p_-,p_+\notin X+w$ and $q_-,q_+\notin X$.  The first assertion follows because subtracting $w$ from either $p_-$ or $p_+$ produces a point with positive imaginary part; the second follows from the choice of the disc around $x_0$ for $q_-$ and from $\Re q_+>M$ for $q_+$.  All four points lie on the indicated boundaries.

For $t>0$, each of the four vertical rays
\[
 p_-+it,\qquad q_-+it,\qquad p_++it,\qquad q_++it
\]
is disjoint from $X\cup(X+w)$.  For the rays from $p_-$ and $p_+$, disjointness from $X$ follows from the supporting-line condition, while disjointness from $X+w$ follows after subtraction of $w$, since the resulting rays lie in the upper half-plane.  For the ray from $q_-$, the part below the supporting line lies in $\overline B(x_0,2y_0)$, while subtracting $w$ gives the ray from $p_-$.  For the ray from $q_+$, the real part is larger than $M$, while subtracting $w$ gives the ray above $r$; the maximal choice of $\Im r$ excludes further points of $X$ on that ray.

Choose a large closed rectangle $R$ containing $X\cup(X+w)$ in its interior and with upper side above the four marked points.  Extend the four vertical rays to that upper side.  By \eqref{eq:access-order}, their endpoints there occur in the cyclic order
\[
        p_-,\ q_-,\ p_+,\ q_+.
\]
Let $T_1,T_2,T_3,T_4$ be the closed $r_0$-neighborhoods of the four access segments, for the ray endpoints $p_-,q_-,p_+,q_+$ in this order; each $T_i$ is convex and contains its marked point in its interior.  Since the access segments lie on distinct vertical lines and are compact and disjoint from the respective opposite sets, $r_0>0$ may be chosen so small that the tubes are pairwise disjoint and
\begin{equation}\label{eq:tube-separation}
        (T_1\cup T_3)\cap(X+w)=\varnothing,
        \qquad
        (T_2\cup T_4)\cap X=\varnothing.
\end{equation}
By compactness, there is $m_0>0$ with $\dist(T_1\cup T_3,X+w)\ge m_0$ and $\dist(T_2\cup T_4,X)\ge m_0$.  Choose $0<d<m_0/4$ so small that the $2d$-neighborhood of $X\cup(X+w)$ lies in $\interior R$.  By \cref{lem:jordan-neighborhood}, choose a compact topological disc $Q$ and $\varepsilon_0>0$ such that
\[
        \mathcal N_{\varepsilon_0}(X)\subset Q\subset\mathcal N_d(X).
\]
Choose $0<\eta_0<d$. Whenever $|w'-w|<\eta_0$, we have
$Q+w'\subset\mathcal N_{2d}(X+w)\subset\interior R$, and the margin $m_0-d-\eta_0>0$ gives
\[
        (T_1\cup T_3)\cap(Q+w')=\varnothing,
        \qquad
        (T_2\cup T_4)\cap Q=\varnothing.
\]
Finally choose open neighborhoods $P_-,V_-,P_+,V_+$ of $p_-,q_-,p_+,q_+$, respectively, with closures contained in the interiors of $T_1,T_2,T_3,T_4$; they are pairwise disjoint because the tubes are.

Let $|w'-w|<\eta_0$, select points $p'_-\in P_-$, $q'_-\in V_-$, $p'_+\in P_+$ and $q'_+\in V_+$, and let $\gamma\subset Q$ join $p'_-$ to $p'_+$ while $\lambda\subset Q+w'$ joins $q'_-$ to $q'_+$.  Extend $\gamma$ to the upper side of $R$ through $T_1\cap R$ and $T_3\cap R$, which are convex sets containing the selected endpoints and meeting the upper side of $R$, and extend $\lambda$ through $T_2\cap R$ and $T_4\cap R$ in the same way.  The endpoints of the extended paths alternate on $\partial R$, so \cref{lem:alternating-endpoints} forces the extended paths to meet.  The tubes are pairwise disjoint; in addition, $T_1\cup T_3$ is disjoint from $Q+w'$ and $T_2\cup T_4$ is disjoint from $Q$.  Hence the intersection must belong to $\gamma\cap\lambda$. This proves all parts of \cref{prop:stable-linked}.

The configuration also satisfies the algebraic intersection condition of \cite[Definition~7.1.3]{CKW2017}. Orient both paths from the minus endpoint to the plus endpoint and compute the intersection number using general-position representatives with the endpoints fixed. This number is unchanged by homotopies in $Q$ and $Q+w'$, since each disc avoids the endpoints of the other path. The alternating boundary endpoints give intersection number $\pm1$ for the extended paths, by comparison with two crossing chords of a disc. The access arcs avoid the opposite paths and contribute no intersections. Thus $\gamma$ and $\lambda$ have nonzero algebraic intersection number.
\end{proof}

\section{Exact verification of the worked trap vector}\label{app:worked-certificate}

This appendix verifies the particular vector used in \cref{ex:worked}.  All terminating decimals in the bounds below denote exact rational numbers.  The argument uses neither a sampled attractor nor a numerical approximation to its filled hull.

\subsection{Root isolation and irrational direction}

Put
\[
        s=\frac{1093336+1366048i}{10^6},\qquad e=10^{-6},
        \qquad P(z)=z^4-2z^3+2z^2+2z-2.
\]
For a Gaussian rational $u$, write $\|u\|_1=|\Re u|+|\Im u|$, so $|u|\le\|u\|_1$.  Direct rational arithmetic gives
\[
 \|P(s)\|_1<3e,\quad |\Re P'(s)|>8,\quad
 \left\|\frac{P''(s)}2\right\|_1<19,\quad
 \left\|\frac{P'''(s)}6\right\|_1<8.
\]
On $|z-s|=e$, the difference between $P(z)$ and $P'(s)(z-s)$ has modulus at most
$\|P(s)\|_1+19e^2+8e^3+e^4<8e<|P'(s)(z-s)|$.
Rouch\'e's theorem therefore gives exactly one zero, counted with multiplicity, in $B(s,e)$.  This zero is simple.  Also $P(-1)=P(1)=1$ and $P(0)=-2$, so $P$ has a real zero in each of $(-1,0)$ and $(0,1)$.  Together with the conjugate disc, these account for all four zeros.  Hence the zero in $B(s,e)$ is the unique root $c_0$ in the upper half-plane.  The same disc gives
\[
        |c_0|>1.749>\sqrt3,\qquad |c_0+1|>\sqrt6,
        \qquad |c_0^{-1}|<0.572.
\]
For the second inequality it suffices to note that
$(\Re s+1-e)^2+(\Im s-e)^2>6$; the other comparisons are likewise rational after squaring positive quantities.

For completeness, the direction of $c_0$ is irrational.  Direct elimination gives
\[
 \operatorname{Res}_z\bigl(P(z),t^4P(z/t)\bigr)=16(t-1)^4R(t),
\]
where
\begin{align*}
 R(t)={}&t^{12}+2t^{11}+4t^{10}+10t^9+37t^8+74t^7+86t^6\\
        &{}+74t^5+37t^4+10t^3+4t^2+2t+1.
\end{align*}
The ratio $t=c_0/\overline{c_0}$ is not $1$ and is a root of $R$.  Modulo $2$, $R(t)=(t+1)^{12}$.  Thus any cyclotomic factor of $R$ must have power-of-two order: if the order is $2^ab$ with $b$ odd, reduction modulo $2$ leaves roots of order $b$, forcing $b=1$.  Its degree is at most $12$, so the only possibilities are
\[
        t-1,\quad t+1,\quad t^2+1,\quad t^4+1,\quad t^8+1.
\]
None divides $R$.  The first three remainders are $342$, $-2$ and $-18$; the last two are
\[
        -62t^3-78t^2-62t
        \quad\text{and}\quad
        74t^7+86t^6+74t^5+36t^4+8t^3-8t-36.
\]
Consequently $c_0/\overline{c_0}$ is not a root of unity, which is equivalent to $\arg(c_0)/\pi\notin\mathbb Q$.

\subsection{A support gap with a rational clearance}

Write $\lambda=c_0^{-1}$ and define
\[
 a_k^*=\begin{cases}2,&\Im(\lambda^k)>0,\\0,&\Im(\lambda^k)<0,\end{cases}
 \qquad
 p=\sum_{k\ge1}a_k^*\lambda^k,\qquad H=\Im p.
\]
There are no zero imaginary parts, by irrationality.  Maximizing the imaginary part digit by digit shows that the upper support line $\ell=\{\Im z=H\}$ meets $K(c_0,3)$ in exactly $p,p+1,p+2$.  It also meets the filled hull $X$ in just these points: $X\subset\conv K(c_0,3)$, and every point of $X$ on $\ell$ lies on $\partial X\subset K(c_0,3)$ by \cref{lem:filled-boundary}.

The following open rational intervals enclose the first eight powers of $\lambda$:
\begin{center}
\begin{tabular}{c@{\qquad}c@{\qquad}c}
$k$ & $\Re(\lambda^k)$ & $\Im(\lambda^k)$\\[2pt]\hline
$1$ & $(0.357,0.358)$   & $(-0.447,-0.446)$\\
$2$ & $(-0.072,-0.071)$ & $(-0.319,-0.318)$\\
$3$ & $(-0.168,-0.167)$ & $(-0.082,-0.081)$\\
$4$ & $(-0.097,-0.096)$ & $(0.045,0.046)$\\
$5$ & $(-0.015,-0.014)$ & $(0.059,0.060)$\\
$6$ & $(0.021,0.022)$   & $(0.027,0.028)$\\
$7$ & $(0.019,0.020)$   & $(0,0.001)$\\
$8$ & $(0.007,0.008)$   & $(-0.009,-0.008)$
\end{tabular}
\end{center}
These bounds follow without root-finding: with $\lambda_*=1/s$, one has $|\lambda-\lambda_*|<e$ and $|\lambda|,|\lambda_*|<0.572$, whence
\[
        |\lambda^k-\lambda_*^k|
        <k(0.572)^{k-1}e.
\]
Substituting the Gaussian rational $\lambda_*$ verifies every endpoint in the table.  The tail bounds
\[
        \frac{2(0.572)^8}{1-0.572}<0.054,
        \qquad
        \frac{2(0.572)^9}{1-0.572}<0.031
\]
will suffice.

Set $x_0=p+1/2$ and $y_0=3/50$.  The signs in the table give
$p=2(\lambda^4+\lambda^5+\lambda^6+\lambda^7)+E$ with $|E|<0.054$, so
\[
        -0.198<\Re p<-0.082,
        \qquad 0.302<\Re x_0<0.418.
\]
Suppose that $z=\sum_{k\ge0}a_k\lambda^k\in K(c_0,3)$ lies in $\overline B(x_0,0.12)$.  The nonnegative support deficits satisfy
\[
        0\le H-\Im z
        =\sum_{k\ge1}(a_k^*-a_k)\Im(\lambda^k)\le0.12.
\]
The first three rows of the table force $a_1=a_2=0$ and $a_3\le1$.  Using the real-part bounds for $3\le k\le8$ and the tail starting at $9$, we obtain
\[
        a_0-0.423<\Re z<a_0+0.131.
\]
Indeed the lower error is at most $0.168+2(0.097)+2(0.015)+0.031=0.423$, and the upper error at most $2(0.022+0.020+0.008)+0.031=0.131$.  Thus $a_0=0$ implies $\Re z<0.131$, while $a_0\ge1$ implies $\Re z>0.577$.  Neither is compatible with the necessary condition $0.182<\Re z<0.538$ for a point of $\overline B(x_0,0.12)$.  This contradiction proves that the closed disc misses $K(c_0,3)$.

The disc is connected and contains points above the support line, hence lies in the unbounded component of $\C\setminus K(c_0,3)$.  It therefore also misses $X$.  Taking $p_-=p$, $p_+=p+1$ and $y_0=3/50$ in \cref{app:linked-translates} gives exactly
\[
        w=x_0-iy_0-p_-=\frac12-\frac3{50}i.
\]
The proof of \cref{prop:stable-linked}, followed by \cref{lem:trap-extraction}, now supplies the required topological disc and open stability margins for this fixed vector.  No claim that the plotted numerical neighborhood equals that disc is needed.

\section*{Acknowledgments}

The author thanks David Juher and Joan Salda\~na for discussions on collinear self-similar sets and for their collaboration on related work, and Christoph Bandt for comments on the neighboring-piece argument.  This work was supported by grant PID2023-146424NB-I00, funded by MICIU/AEI/10.13039/501100011033 and by ERDF/EU, and by the Universitat de Girona and Banco Santander Grant Programme for Researchers in Training (IFUdG 2022--2024).  The author declares no competing interests.

\smallskip
\noindent\textbf{Data availability statement.}
The proofs and the exact inequalities used in this article are contained in the text.  Ancillary files include the figure generators, the numerical mask behind Figure~\ref{fig:regular-core}, and an exact rational-arithmetic checker for \cref{app:worked-certificate}.

\smallskip
\noindent\textbf{AI disclosure statement.}
ChatGPT and Codex (OpenAI) and Claude (Anthropic) were used for editorial revision, development of figure-generation code, and exploration and checking of arguments.  The use of these tools does not replace mathematical verification.  The author is responsible for the arguments, computations, references and final content.

\end{document}